%% file: main.tex
\documentclass[12pt]{elsarticle}

\usepackage{graphicx, import}
\usepackage{amsfonts,textgreek}
\usepackage{mathtools}
\usepackage{dsfont}
\usepackage{multicol}
\usepackage{bbm}
\usepackage{url}
\usepackage[dvipsnames]{xcolor}
\usepackage[colorlinks=true,allcolors=blue]{hyperref}
\usepackage[nameinlink,noabbrev,capitalize]{cleveref}
\usepackage{xfrac}
\usepackage{caption}
\usepackage{subcaption}
\usepackage{xspace}
\usepackage{stmaryrd}
\usepackage{enumerate}
\usepackage{soul}
\usepackage{lscape}

\usepackage[utf8]{inputenc}
\usepackage{wrapfig,bm}
\usepackage{lineno}
\usepackage{siunitx} 
\usepackage{caption}
\usepackage{psfrag}
\usepackage{upgreek}
\usepackage{isomath}
\usepackage{latexsym}
\usepackage{array}
\usepackage{multirow}
\usepackage{booktabs}
\usepackage{colortbl}
\usepackage{verbatim}
\usepackage{graphicx}
\usepackage{epstopdf,overpic}
\usepackage{float}
\usepackage{soul}
\usepackage{color}
\usepackage{dirtytalk}
\usepackage{textcomp}
\usepackage{algorithm, algpseudocode, lmodern, mathtools}
\usepackage{tikz}
\usetikzlibrary{shapes.geometric, arrows}
\usepackage{amssymb}
\usepackage{amsmath}
\usepackage{mathrsfs}  
\usepackage{graphicx}
\usepackage{titlesec}
\usepackage{indentfirst}
\usepackage{bigints}
\usepackage{breqn}
\usepackage[export]{adjustbox}
\usepackage{pgfplots}
\usetikzlibrary{plotmarks,spy,calc}

\usepackage[margin=2.5cm]{geometry}

\usepackage{multirow}

\usepgfplotslibrary{groupplots}

\def\centerarc[#1](#2)(#3:#4:#5);%
    {
    \draw[#1]([shift=(#3:#5)]#2) arc (#3:#4:#5);
    }

\biboptions{sort&compress}

\pgfplotsset{compat=1.15}

\begin{document}
\begin{frontmatter}
\title{Optimal design of Piezoelectric Energy Harvesters for bridge infrastructure: Effects of Location and Traffic Intensity on Energy Production}

\author[label1]{S. Yao}
\author[label1]{P. Peralta-Braz}
\author[label1]{M. M. Alamdari}
\author[label2]{R. O. Ruiz}
\author[label1]{E. Atroshchenko\footnote{Corresponding author, e.atroshchenko@unsw.edu.au, eatroshch@gmail.com}}

\address[label1]{School of Civil and Environmental Engineering, University of New South Wales, Sydney, Australia}
\address[label2]{Department of Mechanical Engineering, University of Michigan-Dearborn, Dearborn, USA}

\journal{Journal of Sound and Vibration}
\begin{abstract}
\input{Sections/abstract}

\end{abstract}


\end{frontmatter}

\section{Introduction}
\label{S:1}
\input{Sections/section1}

\section{Case Study: Cable -Stayed Bridge}
\label{S:2}
\input{Sections/section2}


\section{PEH Modeling}
\label{S:3}
\input{Sections/section3}


\section{Optimisation Framework of PEHs}
\label{3.33}
\input{Sections/section3.1}



\section{Case Study: Optimisation of PEHs for a Cable-stayed Bridge}
\label{S:4}
\input{Sections/section4}

\section{Conclusion}
\label{S:7}
\input{Sections/conclusion}

\section*{Acknowledgements}
This research is undertaken with the assistance of resources and services from the National Computational Infrastructure (NCI), which is supported by the Australian Government. The authors acknowledge support from the University of New South Wales (UNSW) Resource Allocation Scheme managed by Research Technology Services at UNSW Sydney.
\section*{Appendix: }
\setcounter{equation}{0}
\renewcommand\theequation{A.\arabic{equation}}
\addcontentsline{toc}{section}{Appendix:}
\input{Sections/Appendix}
\label{A}

\bibliographystyle{model1-num-names}
\bibliography{reference.bib}

\end{document}

%% file: Sections/abstract.tex
Piezoelectric energy harvesters (PEHs) can be used as an additional power supply for a Structural Health Monitoring (SHM) system. Its design can be optimised for the best performance, however, optimal design depends on the input vibration, e.g. an acceleration of a bridge due to wind loads and passing traffic. In the previous studies, we have shown that an optimal design tunes to some predominant frequency of the input signal spectrum. However, our previous study was limited to a single location of a PEH on a bridge. In this work we extend the rigorous optimisation framework to include the effect of the PEH's location. The optimisation framework uses the PEH Kirchhoff-Love plate model discretised by IsoGeometric Analysis (IGA) and coupled with Particle Swarm Optimisation (PSO) algorithm to find the designs with maximum energy outputs for a large number of input acceleration histories, extracted from the recorded dynamic response data of a real cable-stayed bridge in Australia. Then, clustering is performed to find several best candidates for the entire bridge. Additional study is performed to quantify the effect of traffic intensity on the produced energy. It is shown that variation of traffic volume throughout the 24hr time window leads to variation in optimal PEH designs. The study concludes the impact of location and traffic on energy harvesting by identifying the optimal PEH design and placement throughout the bridge structure. The results indicate that the key factors of maximising energy harvesting efficiency are related to the input excitation and the mode of vibration being excited. The position of maximum displacement in the vibration mode corresponds to the best location for energy harvesting. Also, the best device has a fundamental frequency close to the frequency of the corresponding vibration mode. In addition, the change of traffic intensity affects the amount of convertible mechanical energy and also directs the fundamental frequency of a PEH to shift within a specific range of frequencies to achieve the highest energy conversion.

%% file: Sections/section1.tex
Health monitoring is a necessary component to ensure the normal operation of a wide range of civil and industrial systems. Structural damage and degradation can appear in any location of the structure in various forms and cause serious economic impacts and catastrophic failures, such as the collapse of a bridge. With the popularisation and improvement of monitoring systems, the cost of implementing a structural health monitoring plan and timely prevention is much lower than that of the failure repairment of the entire structure. These factors have promoted the development of Structural Health Monitoring (SHM).

An ideal SHM system should have the ability to detect structural damage immediately and transmit data to locate the fault and analyse the severity so that the corresponding treatment measures (such as stopping operation or repairing) can be taken  \cite{I1}. Although regular manual inspections can obtain relevant data characterising the operational conditions of the system, this time-interval monitoring method is labor-intensive and costly, and it is very likely to have unexpected failures within the monitoring intervals. In addition, some cases are threatening to the safety of inspectors, or even inaccessible, such as inspecting equipment located in the deep sea \cite{27}. The rapid development of electronic systems and ubiquitous use of sensors opened the possibility for continuous health monitoring. The high cost and complicated wiring of wired sensors quickly gave birth to the invention of wireless sensors, and introduced the use of traditional portable batteries \cite{I2}. Although wireless technology eliminates the limitations of the heavy peripheral cable in traditional wired sensors, it brings benefits as well as challenges.

The problems that a perfect health monitoring system should solve include: the number, location, and power supply pattern of the sensors. It is widely accepted that the most critical problem is how to provide continuous and stable energy for the sensor to reduce the cost and labor of frequent replacement of traditional batteries. Therefore, the concept of ''self-powered" have emerged, which refers to the conversion of environmental energy, such as temperature gradient, light and mechanical energy, etc., into useful electrical energy to maintain its own operation and work. Energy harvesting provides a viable and effective alternative to traditional batteries. Common energy harvesting technologies include photovoltaics, electromagnetics, capacitance, pyroelectrics and piezoelectric generation \cite{70}. 

Among them, Piezoelectric Energy Harvester (PEH) that generates electricity based on material deformation caused by external forces is the most suitable choice since mechanical vibrations are ubiquitous \cite{54}, \cite{I3}. Although the design concept of PEH has been studied in depth, its application is still limited due to low conversion efficiency and inability to produce enough power to replace a battery. Therefore, PEH performance optimisation is the key to the exploitation of their ubiquitous use. Benasciutti et al. \cite{126} found that the output power is highly affected by the surface area of the piezoelectric element, especially the width, which proves that shape optimisation is a promising research direction. Roundy et al. \cite{125} proposed a trapezoidal cantilever and experimentally showed that the power output of a beam with a variable width was increased by 30 percent compared to a rectangular beam. Unconventional shapes were studied in \cite{97}, where it was shown that shape optimisation can significantly increase the amount of harvested energy. 


In order to accurately characterize the performance response of a PEH, a reliable and precise mathematical model is needed. In the real world, the excitation received by structures is unpredictable in many cases, such as bridge structures excited by vehicle traffic flow or wind loads. Consequently, many studies need to make assumptions about the sources of excitation \cite{71}. Therefore, most of the performance predictions and experiments on PEHs are based on measurable artificial excitation sources, such as shakers or electrical pulses, rather than random data encountered in real world. These harmonic vibration, which are much simpler than the real excitation, may not be suitable for the optimisation of PEH performance in real operational conditions. Farinholt et al. \cite{a1} designed experiments to simulate vibration input based on the shaker and demonstrated that the powering time for the sensor nodes provided by bridge-based data is significantly different from that provided by harmonic excitation (from 6mins to 30mins). In order to further improve the energy output, Paul et al. \cite{a2} experimentally designed a PEH tuned to one of the bridge's natural frequency and obtained a desired output. 

While superior energy collection capabilities were achieved, it was shown that solving non-convex problem to explore optimal designs was not feasible due to the limitation of fabrication. Therefore, numerical models that are more flexible than experiments have been developed to perform behavioral simulations and performance evaluations of PEHs. However, the differences in multiple theoretical models and analysis methods of PEHs interfere with the accuracy and validity of the simulation results. In order to simulate the output performance of PEHs, many previous studies modeled the energy harvesting system as a lumped parameter model (Single Degree of Freedom model) that couples second-order differential equations with lumped parameters and circuit equations through piezo-electric constitutive equation, which produces inaccurate results \cite{75}. With the expansion of related research, some correction coefficients \cite{77} and models with distributed mass based on beam theories \cite{79,81} were proposed to improve the shortcomings of lumped parameter models, but in fact, the geometric shape defined by the low aspect ratio of cantilever PEHs is more similar to a plate than a beam. Therefore, the theoretical PEH model based on Kirchhoff-Love plate and Isogeometric Analysis (IGA) was proposed in \cite{97} to provide parametric descriptions of complex shapes and accurate solutions. The accuracy and efficiency of the IGA PEH model are attributed to the advantages of IGA with high-order Non-Uniform Rational B-Splines (NURBS) to parameterise the computational domain exactly (thus, avoiding domain discretisation error) and represent higher-order derivatives appearing in the weak form of the Kirchhoff-Love plate model. Thus, higher accuracy is achieved with lower computational burden. The feasibility of the IGA PEH model to replace the experiment was also verified by comparing the dynamic responses of simulated and experimental PEHs under the same excitation conditions.  An extended study of shape optimisation for devices with variable thicknesses was further performed in \cite{extra1}, demonstrating the broad applicability of the PEH IGA model. 

In the subsequent work, \cite{ex35}, the model was used to propose a shape optimisation framework for the actual operating environment, based on the structural health monitoring data of a cable-stayed bridge in New South Wales, Australia. The framework performs shape optimisation for a large number of the so-called {\it events}, i.e. 30 seconds acceleration histories due to passing traffic or wind loads. Then clustering the optimal geometries to identify the best candidates for continuous energy generation.     


However, the above study was limited to a single fixed location of the PEH on the bridge, which is not the best option for maximising the power generation of the entire structure since the sensors are distributed at different locations on the bridge. Hence, in order to maximise the energy harvesting efficiency of PEH systems in practical infrastructure applications, it is important to understand the influence of the PEH location on the harvesting capacity of PEH and identify locations with maximum harvesting capacity. 

Some other studies addressed these questions. For example, Zang et al. \cite{ex40} fabricated an experimental platform including vehicle-bridge interactions to simulate the energy harvesting potential of two PEHs (based on the natural frequencies of the bridge and vehicle-bridge coupling vibration, respectively) under real operating conditions, and the results showed that the device located in the middle of the bridge collect higher energy. But the studies in \cite{ex40} only considered one condition, i.e., a single location point or a single collector configuration, which failed to maximise the energy harvesting efficiency of PEH systems in practical applications. Karimi et al.\cite{4.3} studied the power spectral density diagrams of the system containing the car and the bridge, and found that the power generation of PEHs had two peaks at the natural frequency of the structure. Song et al.\cite{ex41} verified the above statement by simulating the vertical acceleration response induced by a moving train at different locations of two bridges (single-span and double-span).

The objective of this work is to extend the PEH shape optimisation framework proposed in \cite{ex35}, to multiple locations across the bridge. The dynamic response of a real cable-stayed bridge in Australia \cite{4.2} is used to identify optimal PEHs' locations on the bridge and study the impact of PEHs' locations on their energy harvesting performance. The proposed framework provides a benchmark for the future design and installation of PEHs on a bridge. Moreover, it can be seen from the optimisation results of Peralta et al. \cite{ex35} that the optimal shape designs corresponding to different traffic events exhibit high dispersion. But in real life, we tend to design one or a few  devices which would be optimal for multiple locations and operating conditions. Hence, another objective of this work is to propose a design framework that allows to cluster PEHs' locations, thereby avoiding a plethora of optimal designs for ease of installation and fabrication while maximising energy output. Additionally, the work explores the impact of traffic on optimal locations and designs by considering variations of traffic within 24 hours.

The objectives of this work are formulated as follows:
\begin{itemize}
    \item Explore the impact of locations on power generation and find the optimal locations based on the data collected from a real bridge
    \item Design the best PEH at each location and cluster them to reduce the number of optimal designs for the entire bridge
    \item Analyse the impact of traffic on energy harvesting and optimal designs over 24 hrs time windows
\end{itemize}

The remainder of this paper is organised as follows. The cable-stayed bridge used to provide real data is presented in Section \ref{S:2}. The theoretical model of PEH based on the isogeometric analysis applied in this work is explained in Section \ref{S:3}. The optimisation framework of PEHs is presented in Section \ref{3.33}. Next, Section \ref{S:4} presents the results of applying the optimisation procedure and analyses the impact of location and traffic on energy harvesting. Conclusions are presented in Section \ref{S:7}.

%% file: Sections/section2.tex
\begin{figure}[H]
	\centering
	\includegraphics[scale=0.7]{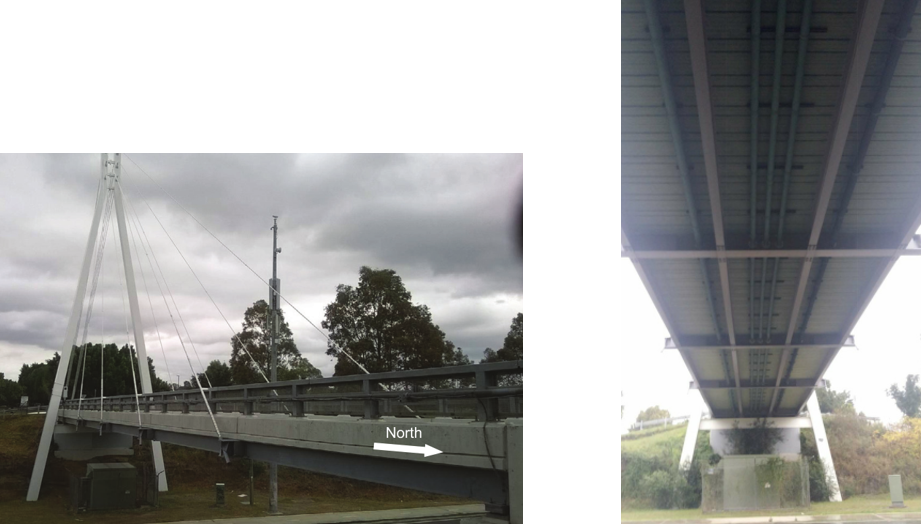}
	\caption{A cable-stayed bridge over the Great Western Highway in the state of New South Wales in Australia \cite{kalhori}.}
	\label{case}
\end{figure}

The case study in this paper is based on an in-service cable-stayed bridge (Figure \ref{case}) consisting of one traffic lane and one pedestrian lane with a maximum load of 30t located on the Great Western Highway in the state of New South Wales, Australia \cite{4.2}. The bridge uses sixteen steel stay cables (38mm in diameter) fixing the A-shaped steel pylon to form an asymmetric bridge structure. The bridge deck (composite steel-concrete deck) with a width of 6.30m, a thickness of 0.16m, and a span of more than 46m is supported by four longitudinal girders (LG1:LG4) and seven cross girders (CG1:CG7) connecting them.

\begin{figure}[H]
	\centering
	\includegraphics[width=\linewidth]{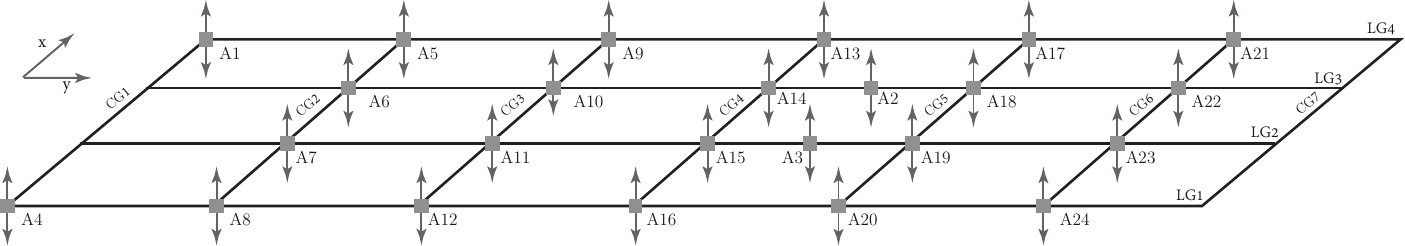}
	\caption{Schematic diagram of the distribution of accelerometers \cite{4.2}.}
	\label{4.5}
\end{figure}

In 2016, an array of uniaxial accelerometer sensors, placed at different positions of the bridge deck as shown in Figure \ref{4.5}, was employed to continuously measure the acceleration response data. The information is sampled at 600Hz and transmitted to the database via a 4G cellular network. The distribution positions of the 24 sensors are shown in Figure \ref{4.5}. Their specifications are all low-noise Silicon Designs accelerometers (Model Number 2210-002) with a sensitivity of 2000 mV/g. This database was used in \cite{ex35} for PEH design optimisation at fixed sensor location A14. In the present study, we extend the design optimisation framework to all sensor locations: A5 - A20. Since the devices located at fixed supports (A1, A4, A21, A22, A23, and A24) will produce a negligible amount of energy, these locations are excluded from consideration.

\begin{figure}[H]
	\centering
	\includegraphics[width=0.85\linewidth]{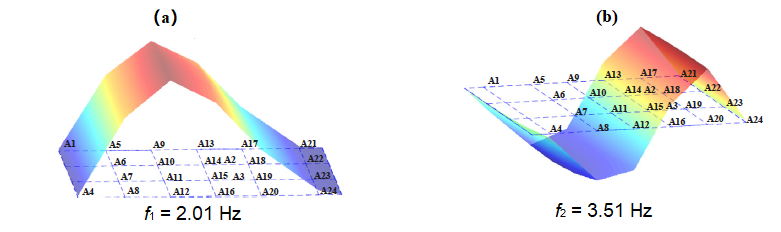}
	\caption{Illustration of the first two bending modes corresponding to the first two natural frequencies obtained from field test: (a) - first mode, (b) - second mode \cite{ex37}.}
	\label{F4.2.3}
\end{figure}

In what follows, we will also need to analyse first two bending modes of the bridge, corresponding to first two natural frequencies: $f_1 = 2.01$ Hz and $f_2 = 3.51$ Hz. They are shown in Figure \ref{F4.2.3} together with the network of accelerometers.


%% file: Sections/section3.tex
\begin{figure}[h]
	\centering
	\includegraphics[width=1\linewidth]{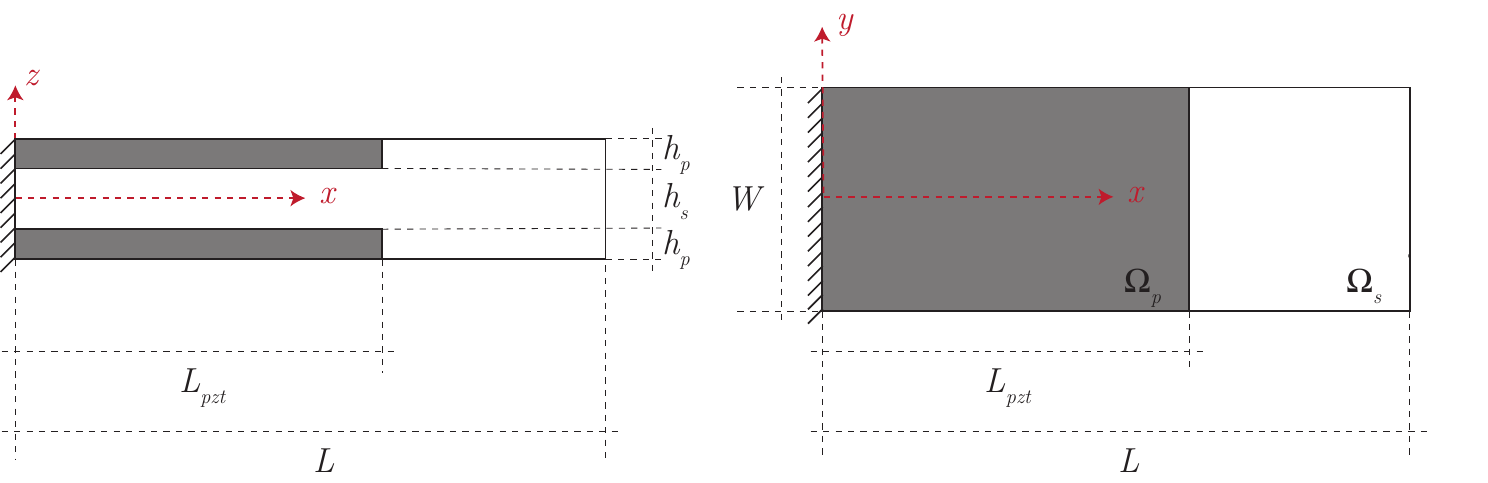}
	\caption{Schematic diagram of a PEH.}
	\label{4.1}
\end{figure}

In this section, we briefly recall the PEH IGA model, originally developed in \cite{97} and further equipped with Model Order Reduction (MOR) and time integration algorithms in \cite{ex35}. The PEH is modeled as a three-layer cantilevered Kirchhoff-Love plate, consisting of two layers of piezoelectric ceramics connected in series to an external electrical resistance $R_{l}$ and a non-piezoelectric structure between them to increase stiffness, as shown in Figure \ref{4.1}. Five parameters define the geometric configuration of the PEH: the total length $L$, the total width $W$, the length $L_{pzt}$ and thickness $h_{p}$ of the piezoelectric layer, and the thickness $h_{s}$ of the structural layer. The thickness of all layers is constant with the total thickness being  $h = 2h_{p} + h_{s}$. In addition, the piezoelectric layer $\Omega_{p}$ partially covers the non-piezoelectric layer $\Omega_{s}$, which was shown to yield better energy production in comparison with a fully covered model \cite{97}. The width $W$ is assumed constant. Although the addition of a tip mass is beneficial for improving the excitation amplitude in the low-frequency region, since this work focuses on the change in the distribution mass caused by the geometric change of PEH, the effect of tip mass is not considered ($M_{tip}$ = 0).

\subsection{IGA Model of a PEH}
The theoretical model to simulate PEH behaviour applied in this work is based on the Kirchhoff-Love plate and the generalised Hamilton’s principle, and its numerical solution is obtained by Isogeometric Analysis \cite{97}. The constitutive relation of the non-piezoelectric substructure layer, assumed to be composed of isotropic homogeneous material under the action of plane-stress, is given by
\begin{equation}
    \mathbf{T} = \mathbf{c}_s \mathbf{S}
    \label{e1}
\end{equation}
where $\mathbf{T}$, $\mathbf{S}$, and $\mathbf{c}_s$ are the mechanical stress, mechanical displacement, and elastic stiffness matrix, respectively. In addition, the constitutive equation of the piezoelectric material is as follows,
\begin{equation}
    \begin{bmatrix}
    \mathbf{T} \\
    \mathbf{D}
    \end{bmatrix}
    =
    \begin{bmatrix}
    \mathbf{c}_{p}^E    &   -\mathbf{e}^T \\
    \mathbf{e}          & \varepsilon^S 
    \end{bmatrix}
    \begin{bmatrix}
    \mathbf{S} \\
    \mathbf{E}
    \end{bmatrix}
    \label{e2}
\end{equation}
where $\mathbf{D}$ represents the electrical displacement, $\mathbf{E}$ is the electric field, $\mathbf{c}_{p}^E$, $\varepsilon$, and $\mathbf{e}$ denote the elastic stiffness matrix under constant electric field, permittivity components matrix at constant strain, and piezoelectric constant matrix, respectively. 

Non-Uniform Rational B-Splines (NURBS) basis functions $N_I$ are used to represent the shape $\mathbf{x}$ and approximate the vertical displacement $w$ in parametric space $\xi\in[0,1]\times[0,1]$: 
\begin{equation}
    \mathbf{x} (\xi) = \sum_{I=1}^k N_I(\xi) \mathbf{\widetilde{x}_I}
    \label{e3}
\end{equation}

\begin{equation}
    {w} (\xi,t) = \sum_{I=1}^k N_I(\xi) {w_I(t)}
    \label{e4}
\end{equation}
where $I$ represents the global index of control points $\mathbf{\widetilde{x}_I}$ and $w_I(t)$ defines the deflection projected at control point $I$. Then by coupling Equation \ref{e3}, Equation \ref{e4}, and the weak form of the problem, a system of algebraic equations consisting of the mechanical equation and the electric equation is derived.
\begin{equation}
    \mathbf{M\ddot{w}}+\mathbf{C\dot{w}}+\mathbf{Kw}-\mathbf{\Theta} v_p=\mathbf{F}a_b
    \label{e5}
\end{equation}
\begin{equation}
    C_p\dot{v}_p+\frac{v_p}{{R}_l}+\mathbf{\Theta^T \dot{w}}=0
    \label{e6}
\end{equation}
where $\mathbf{w} \in \mathbb{R}^{N\times 1}$, $\mathbf{M} \in \mathbb{R}^{N\times N}$, $\mathbf{K} \in \mathbb{R}^{N\times N}$, $\mathbf{C}=\alpha \mathbf{M}+\beta \mathbf{K}\in \mathbb{R}^{N\times 1}$, $\mathbf{F} \in \mathbb{R}^{N\times 1}$ and $\mathbf{\Theta} \in \mathbb{R}^{N\times 1}$  denote the vector of deflections $w_I(t)$, the mass matrix, the stiffness matrix, the mechanical damping matrix, the mechanical forces vector, the electromechanical coupling vector, respectively. $\alpha$ and $\beta$ represent the proportional damping coefficients. $a_b$ defines the input acceleration, and $v_p$ represents the output voltage. $C_p$ and ${R}_l$ are the capacitance and the external
electric resistance, respectively. The full expressions of matrices $\mathbf{M}$, $\mathbf{K}$, $\mathbf{\Theta}$, and $\mathbf{F}$ are given in the Appendix.  

Based on Equations (\ref{e5}) and (\ref{e6}), the Frequency Response Function (FRF) shown below can be derived from the relationship between output voltage ($v_p=V_o e^{i\omega t}$) and input acceleration ($a_b=A_b e^{i\omega t}$) to analyse the dynamic behavior of PEH during energy harvesting under harmonic excitations, where $i=\sqrt{-1}$ \cite{ex35}.
\begin{equation}
    \begin{aligned}
    H_v(\omega) & = \frac{V_o(\omega)}{A_b(\omega)} \\ &
    = i\omega \left( \frac{1}{R_l} + i\omega C_p \right)^{-1} \mathbf{\Theta}^T \left(-\omega^2 \mathbf{M} + i\omega \mathbf{C} + \mathbf{K} +i\omega \left(\frac{1}{R_l} +i\omega C_p \right)^{-1} \mathbf{\Theta} \mathbf{\Theta}^T \right)^{-1} \mathbf{F}
    \end{aligned}
    \label{e22}
\end{equation}

\subsection{Modal Order Reduction (MOR) and Time Integration}
Large-scale discretisation in simulating the dynamic behavior of PEH greatly increases the computational load and time. A feasible solution is Modal Order Reduction (MOR)\cite{ex34}, which can effectively reduce the number of degrees of freedom. A truncated expansion of the $K$ first modal vectors is applied to approximate the deflection solution $\mathbf{w}_o$ as 
\begin{equation}
    \mathbf{w}_o = \mathbf{\Phi}_o\bm{\eta} 
    \label{e14}
\end{equation}
where $\mathbf{w} _o \in \mathbb{R}^{N\times1}$ represents the approximate deflection in terms of modal coordinates $\bm{\eta}\in\mathbb{R}^{K\times1}$. And the $K$ first mode shape vectors $\phi _i$ constitute matrix $\mathbf{\Phi}_o\in \mathbb{R}^{N\times K}$.

Then, the reduced order system can be obtained from Equations (\ref{e5}) and (\ref{e6}) normalised by the reduced mass matrix as
\begin{equation}
    \mathbf{\ddot{\bm{\eta}} } +\mathbf{c}_o\mathbf{\dot{\bm{\eta}}}+\mathbf{k}_o\mathbf{\bm{\eta}} -\bm{\theta} _ov_p = \mathbf{f}_o a_b
    \label{e16}
\end{equation}

\begin{equation}
    C_p\dot{v_p}+\frac{v_p}{R_l}+\mathbf{\Theta} ^T\mathbf{\Phi} _o\mathbf{\dot{\bm{\eta}}}=0
    \label{e17}
\end{equation}

where

\begin{equation}
    \mathbf{c}_o=\mathbf{m}_o^{-1}\mathbf{\Phi} _o^T\mathbf{C}\mathbf{\Phi} _o=\begin{bmatrix}
2\zeta _1w_1 &  & \\ 
 & \ddots  & \\ 
 &  &2\zeta _K w_K 
\end{bmatrix},
\quad
\mathbf{k}_o=\mathbf{m}_o^{-1}\mathbf{\Phi} _o^T\mathbf{K}\mathbf{\Phi} _o=\begin{bmatrix}
w_1^2 &  & \\ 
 & \ddots  & \\ 
 &  & w_K^2 
\end{bmatrix}\label{e18}
\end{equation}


\begin{equation}
    \bm{\theta} _o = \mathbf{m}_o^{-1}\mathbf{\Phi} _o^T\bm{\theta}, \quad \mathbf{f} _o = \mathbf{m}_o^{-1}\mathbf{\Phi} _o^T\mathbf{F}, \quad \mathbf{m}_o = \mathbf{\Phi}_o^T \mathbf{M} \mathbf{\Phi}_o 
    \label{e20}
\end{equation}


The solution of Equations (\ref{e16}) and (\ref{e17}) represent the voltage output of the PEH. Time integration of Equations (\ref{e16}) and (\ref{e17}) is performed using Simulink ode45 solver. This solver uses the Runge-Kutta (RK) method \cite{RK}. In \cite{ex35}, its accuracy was demonstrated by performing voltage measurement experiments based on a commercial PEH. The piezoelectric layers of the PEH are connected in series to a 1000 $\Omega $ electrical resistance. The measured dynamic response based on pulse excitation was shown to be very consistent with the numerical simulation results after the modal order reduction, which proves the validity and feasibility of the proposed model.

%% file: Sections/section3.1.tex

\begin{figure}[H]
	\centering
	\includegraphics[scale=0.42]{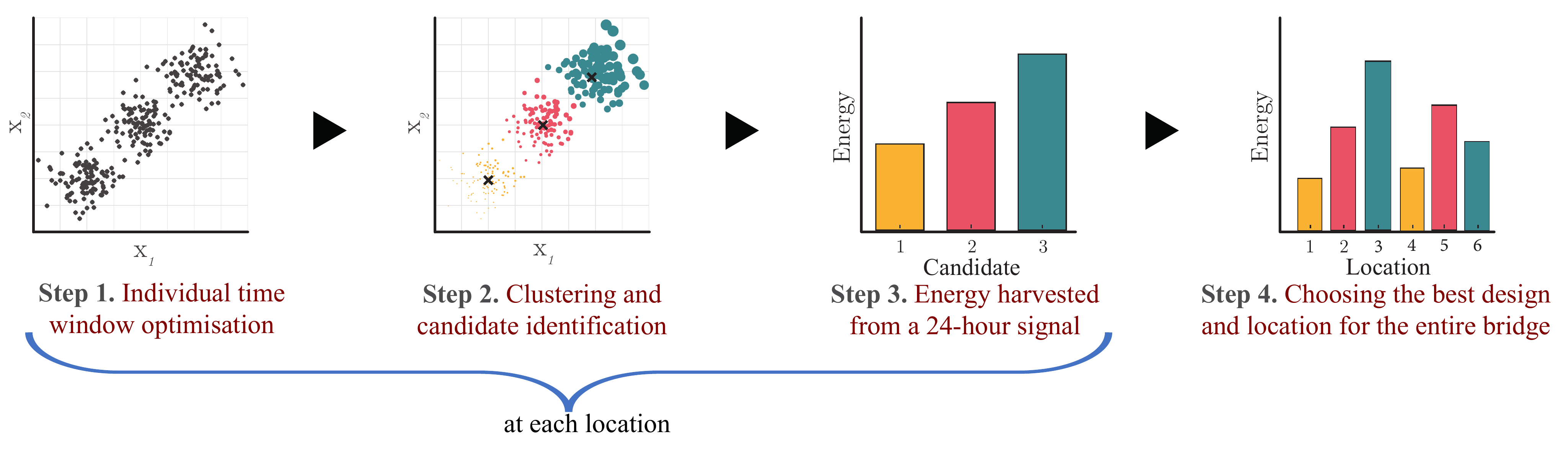}
	\caption{Optimisation framework.}
	\label{step}
\end{figure}

The study adopts the optimisation framework proposed by Peralta et al \cite{ex35} to investigate the maximum potential of a PEH and the effect of its position on a bridge. The framework is divided into four parts, as shown in Figure \ref{step}. The first three steps are carried out at different bridge locations, and the results are compared in the final step. The aim of \textbf{Step 1} is to determine the optimal geometry of PEH $\bold{x^*}$ for various acceleration time windows, which are denoted as $e$ and have a duration of one hour ($T$= 3600sec), representing different times in a single day. The optimisation problem can be expressed as follows:
\begin{equation}
    \bold{x^*} = \arg \underset{\textbf{x}\in \textbf{X}}{max}\left | E(\bold{x}\mid{e})\right |
    \label{e5.1}
\end{equation}
where $E(\bold{x}\mid e)$ is the energy for a PEH configuration given by the geometrical parameters $\bold{x} \in\bold{X}$ harvested from time window $e$, defined as the integration of squared voltage $v_p = v(t,\bold{x}\mid e)$ scaled by numerical resistance $R_l$:
\begin{equation}
    E(\bold{x}\mid e) = \int_{0}^{T}\frac{v^2( t,\bold{x}\mid e)}{R_l}dt
    \label{eE}
\end{equation}

In \textbf{Step 2}, the optimal PEH configurations obtained in the previous step are grouped based on their geometrical parameters (design variables) to reduce the number of candidates. The $k-$means method \cite{ex38} is used, where the centroid of each cluster is calculated by minimising the distance between each data point and the centroid. The number of groups, $k$, is determined using the Silhouette Coefficient method \cite{ex39}. 

The objective of \textbf{Step 3} is to select the optimal geometry for continuous energy generation. In a practical setting, it is desirable to have a device that can produce the maximum energy without being affected by bridge traffic. Hence, the optimal design is chosen from the top $k$ candidates obtained in the previous step based on the energy produced during the 24-hour continuous acceleration signal.


Finally, in \textbf{Step 4}, the impact of location on energy production for the entire bridge is analysed. The optimal geometries from each location obtained in the previous step are compared and ranked, thereby determining the specific zone with the highest energy generation.

%% file: Sections/section4.tex
In this section, we apply the proposed optimisation framework to the cable-stayed bridge, introduced in Section \ref{S:2}. 
As mentioned in Section \ref{S:3}, this paper adopts a bimorph three-layer PEH cantilever plate, the inner non-piezoelectric support layer is made of bronze, and the upper and lower piezoelectric layers are made of PZT 5A ceramics. All material properties and parameters are given in Table \ref{T4.1} and Table \ref{T4.2}. Four parameters are chosen to represent the design space of the PEH shape: the device length $L$, the thickness ratio $H = h_p/h$, the length ratio $l = L_{pzt}/L$, and the aspect ratio $R = W/L$. The range of the above variables is set to $L \in [0.1, 0.5]$ m, $H \in [0.05, 0.45]$, $l \in [0.1, 1]$, $h = 1$ mm. Studies in \cite{ex35} demonstrated that optimal designs correspond to $R = 1$, which represents a square shape. This is due to the fact that the change of $R$ hardly affects the fundamental frequency of the PEH device, while the increase of $R$ and $L$ can promote higher conversion of energy due to the increase of area and distributed mass \cite{97}. Hence, for this study we set $R = 1$, and the shape of a PEH is therefore parameterised in the space of three variables ($L - l - H$), such that the design space is defined as $\bold{x}=[L\,\,l\,\,H]$.


\subsection{Shape optimisation at one location} \label{4.2.1}

\begin{figure}[H]
	\centering
	\includegraphics[width=1\linewidth]{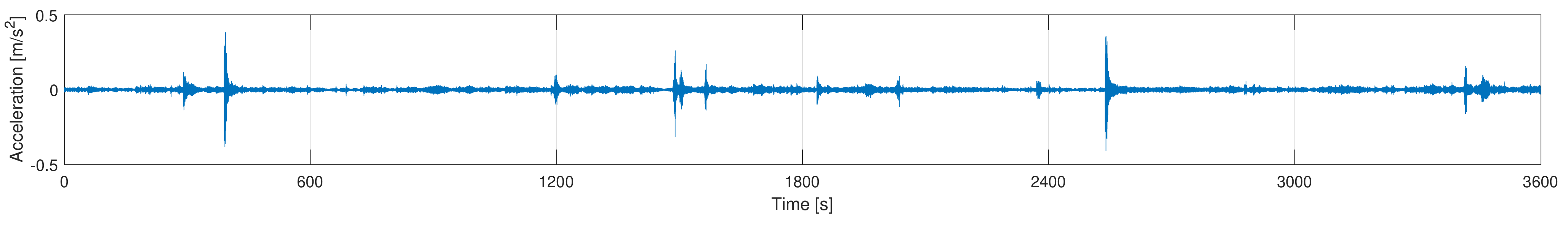}
	\caption{Continuous time window with one-hour acceleration signal observed at sensor location A2.}
	\label{5.8}
\end{figure}

\begin{figure}[H]
	\centering
	\includegraphics[width=1\linewidth]{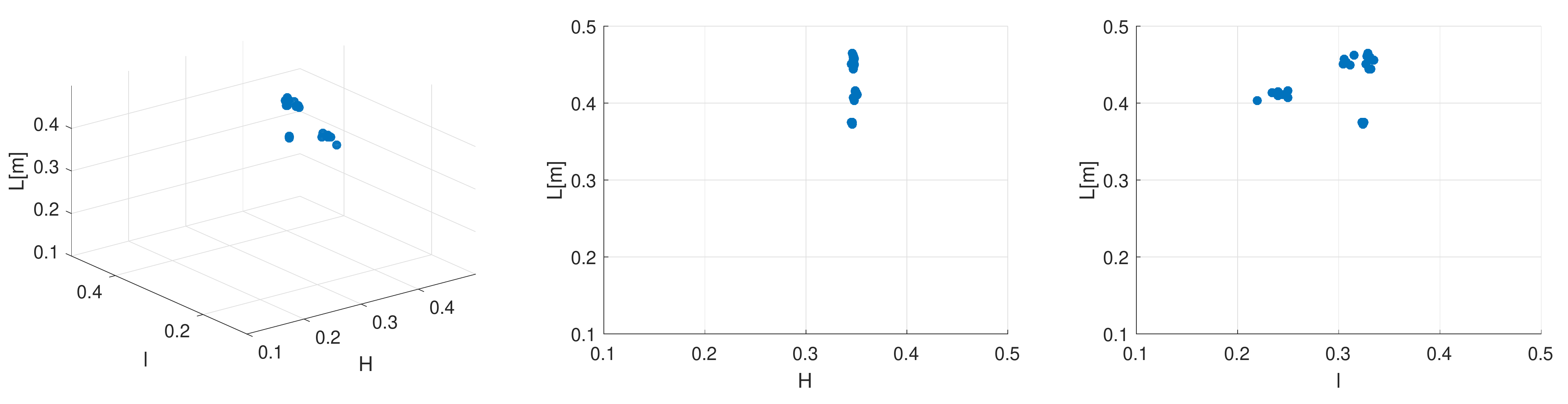}
	\caption{Optimal designs at sensor location A2 for 24 windows with 1-hour duration in the $L - l - H$ space (left), $L - H$ plane (middle) and $L - l$ plane (right).}
	\label{4.1.1}
\end{figure}

The optimisation framework outlined in section \ref{3.33} is applied to all sensor locations on the bridge structure, ranging from A5 to A20. In this demonstration, location A2 is selected for analysis. Figure \ref{5.8} presents the acceleration response for a typical one-hour window extracted from the continuous bridge acceleration observed at sensor location A2. The optimal PEH shapes at A2 are depicted in Figure \ref{4.1.1} using three-parameter spaces: full ($L - l - H$) space and its projections on the $L - l$ and $L - H$ planes. The analysis is performed for 24 such time windows. The results presented in Figure \ref{4.1.1} show some variation in the optimal designs for the parameters $L$ and $l$, while relatively low variation is observed for parameter $H$. The discussion of this observation will be further resumed.

\begin{table}[h]
    \centering
    \renewcommand\arraystretch{1.5}
    \caption{General characteristics of the PEH extracted from \cite{4.1}.}
    \begin{tabular}{p{0.75\linewidth} p{0.15\linewidth}}
    \hline
        Length, $L$ [mm]                              &  $L$ \\
        Dimensionless piezoelectric length, $l$                             &  $l$ \\ 
        Dimensionless piezoelectric thickness, $H$                             &  $H$ \\
        Aspect ratio, $R$                             &  1 \\
        Width, $W$ [mm]                               &  $R L$ \\
        Total thickness of the device, $h$ [mm]       &  1 \\
        Thickness of the PZT layers, $h_p$ [mm]       &  $H h$ \\
        Length of the PZT layers, $L_{pzt}$ [mm]       &  $l L$ \\ 
        Thickness of the substructure layers, $h_s$ {[}mm{]} & $h - 2h_p$ \\ 
        Substrate Young's Modulus, $E$ [GPa]          &  105  \\
        Substrate Density, $\rho_s$ [kg/m$^3$]        &  9000 \\
        Substrate Poisson's ratio, $\nu$              &  0.30 \\
        PZT-5A Density, $\rho_p$ [kg/m$^3$]           &  7800 \\
        PZT Poisson's ration, $\nu$                   &  0.30 \\
        Damping Coefficient, $\alpha$ [s/rad]         &  14.65 \\
        Damping Coefficient, $\beta$ [s/rad]          &  $1\times 10^{-5}$ \\ \hline
    \end{tabular}
    \label{T4.1}
\end{table}

\begin{table}[h]
    \centering
    \renewcommand\arraystretch{1.5}
    \caption{Electro-mechanical properties of the piezoelectric layers extracted from \cite{85}.}
    \begin{tabular}{p{0.75\linewidth} p{0.15\linewidth}}
    \hline
        $c_{11}^{E}, c_{22}^{E}$ [GPa]                &  69.5 \\
        $c_{12}^{E}$ [GPa]                            &  24.3 \\
        $c_{66}^{E}$ [GPa]                            &  22.6 \\
        $d_{31}^{E}$ [C/m$^2$]                        &   -16 \\
        $d_{32}^{E}$ [C/m$^2$]                        &  -16 \\
        Permittivity, $\epsilon_{33}$ [nF/m]          &  9.57 \\
    \hline
    \end{tabular}
    \label{T4.2}
\end{table}

\begin{figure}[H]
	\centering
	\includegraphics[width=0.7\linewidth]{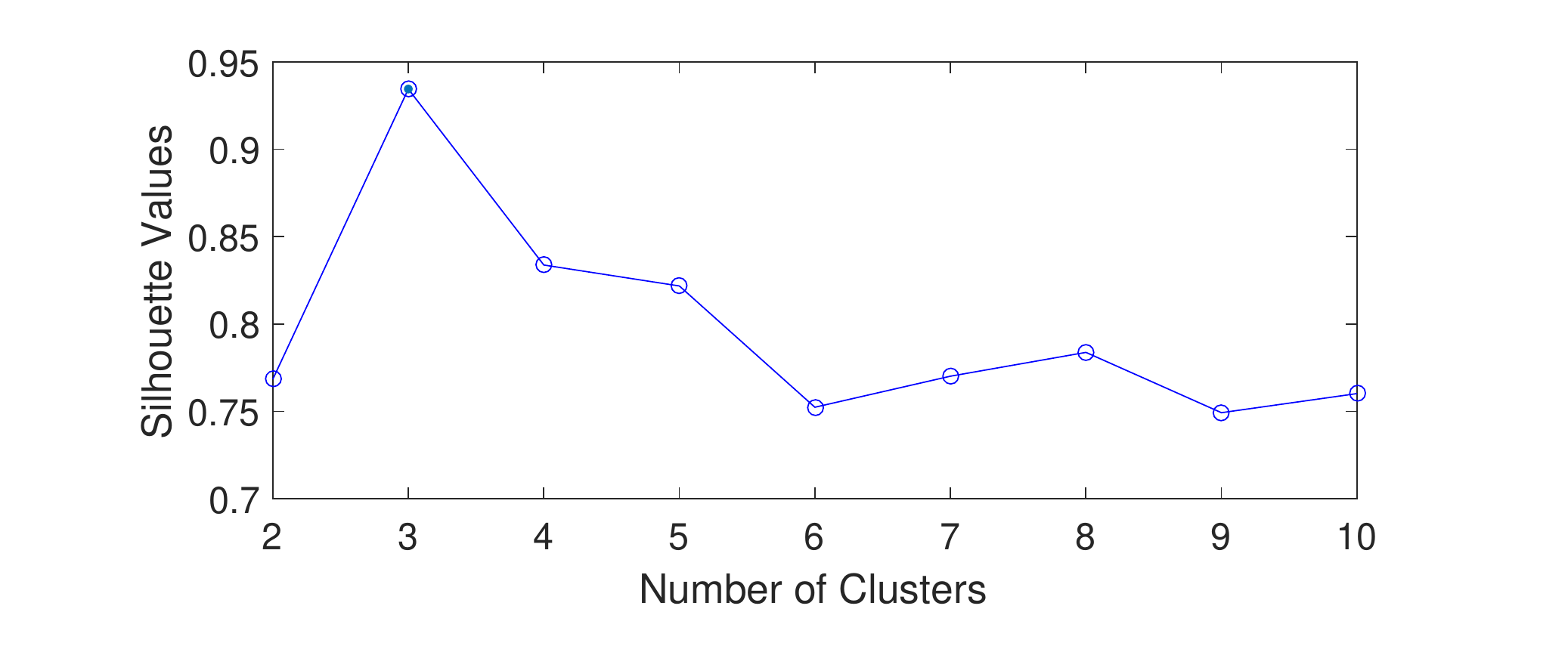}
	\caption{Silhouette Coefficient for the optimal geometries at sensor location A2.}
	\label{4.1.2}
\end{figure}

\begin{figure}[H]
	\centering
	\includegraphics[width=1\linewidth]{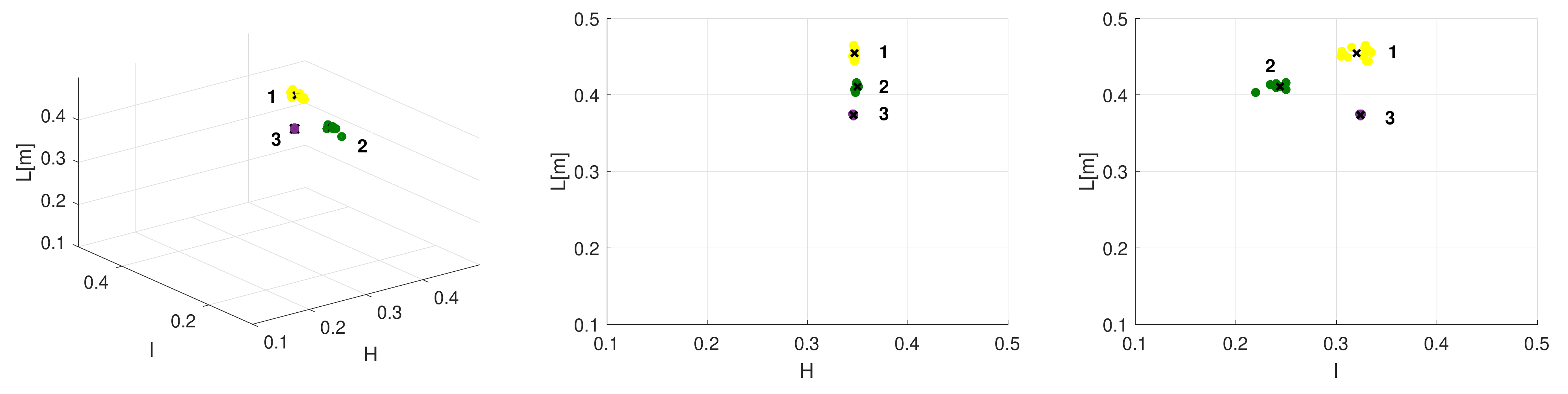}
	\caption{Three clusters of optimal PEH designs.}
	\label{4.1.3}
\end{figure}

The next step in the framework is to simplify the set of candidate designs by grouping the optimal PEHs through $k$-means clustering. This reduces the number of possibilities to a smaller set of representative geometries, with the cluster centroids serving as the best candidates. To determine the number of clusters, the Silhouette Coefficient is evaluated, as seen in Figure \ref{4.1.2}. The highest Silhouette Value determines the optimal number of clusters, which was found to be $k=3$. In Figure \ref{4.1.3}, sample points belonging to the three clusters are indicated by yellow (cluster 1), green (cluster 2), and purple (cluster 3) colors, while the cluster centroids are marked with crosses.


\begin{figure}[H]
	\centering
	\includegraphics[width=0.65\linewidth]{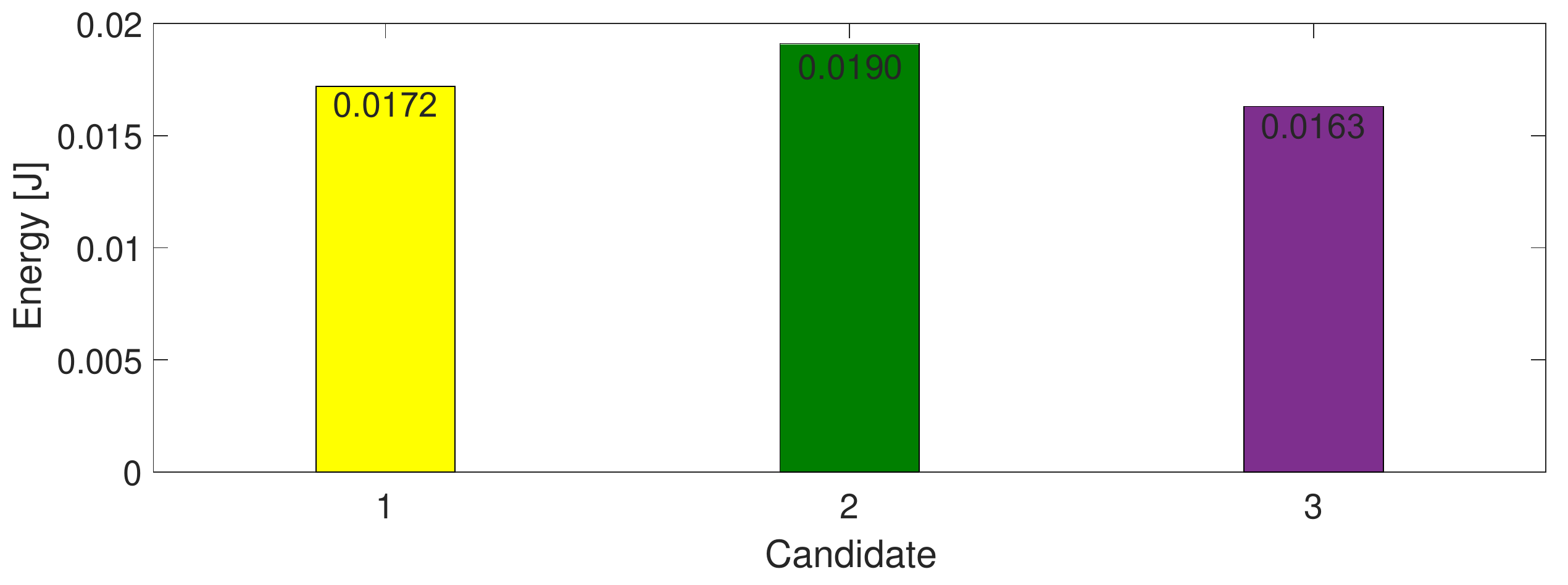}
	\caption{Energy harvested from a 24-hour time window at sensor location A2 by three best candidates.}
	\label{4.1.5}
\end{figure}

The next step involves the evaluation of the energy harvesting of the three best PEH designs based on the 24-hour continuous excitation signal extracted from the recorded acceleration data. Figure \ref{4.1.5} displays the 24-hour energy outputs of the three candidates at sensor location A2. The results show that the PEH design from cluster 2 collects the most energy, with a 10.5\% and 16.6\% increase compared to candidates 1 and 3, respectively. This demonstrates that candidate 2 is the optimal PEH design for the sensor location A2. 

\begin{figure}[ht!]
    \centering
    \begin{subfigure}[b]{0.28\textwidth}
        \includegraphics[width=\textwidth]{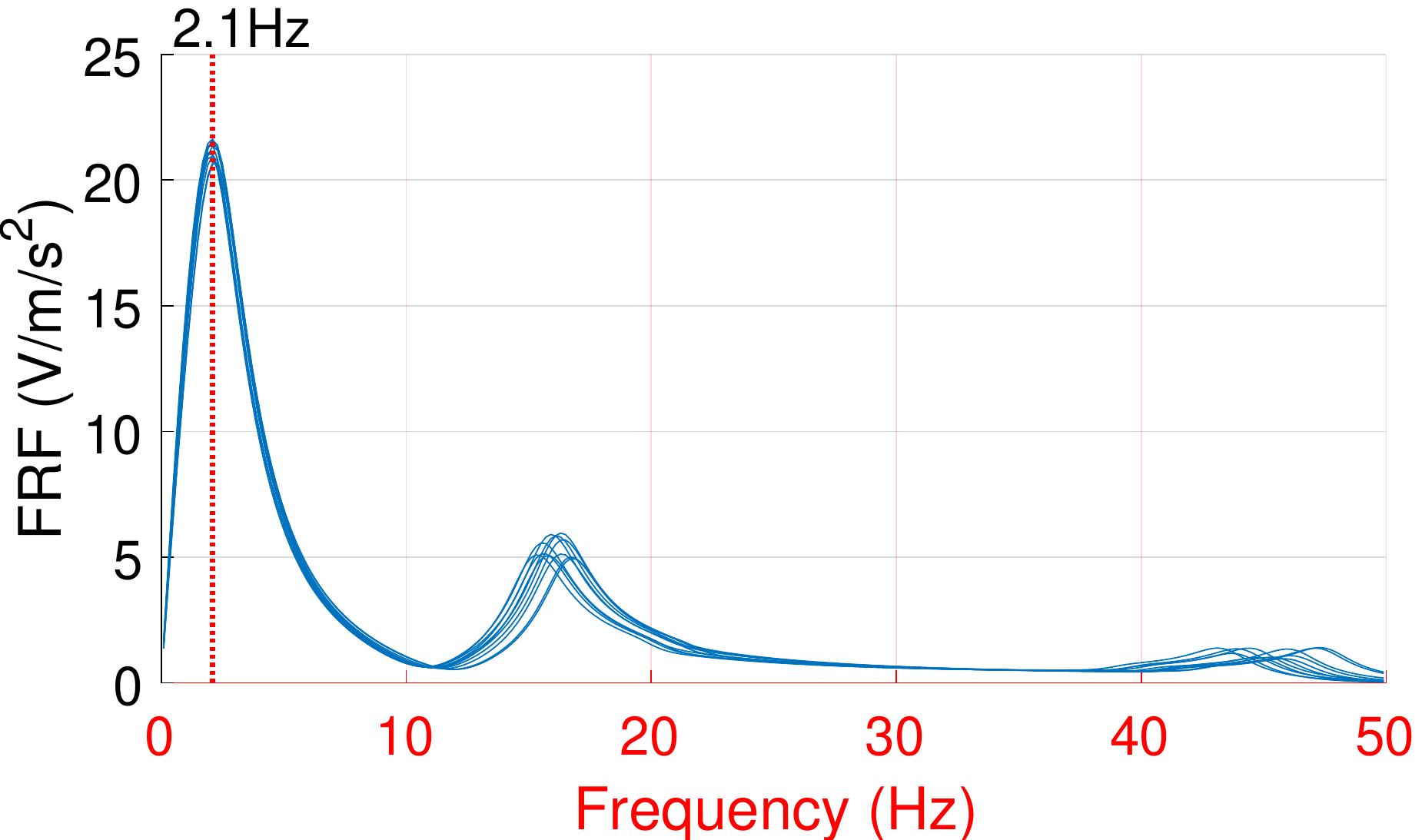}
        \caption{FRFs of all PEHs in cluster 1.}
    \end{subfigure}
    \quad
    \begin{subfigure}[b]{0.28\textwidth}
        \includegraphics[width=\textwidth]{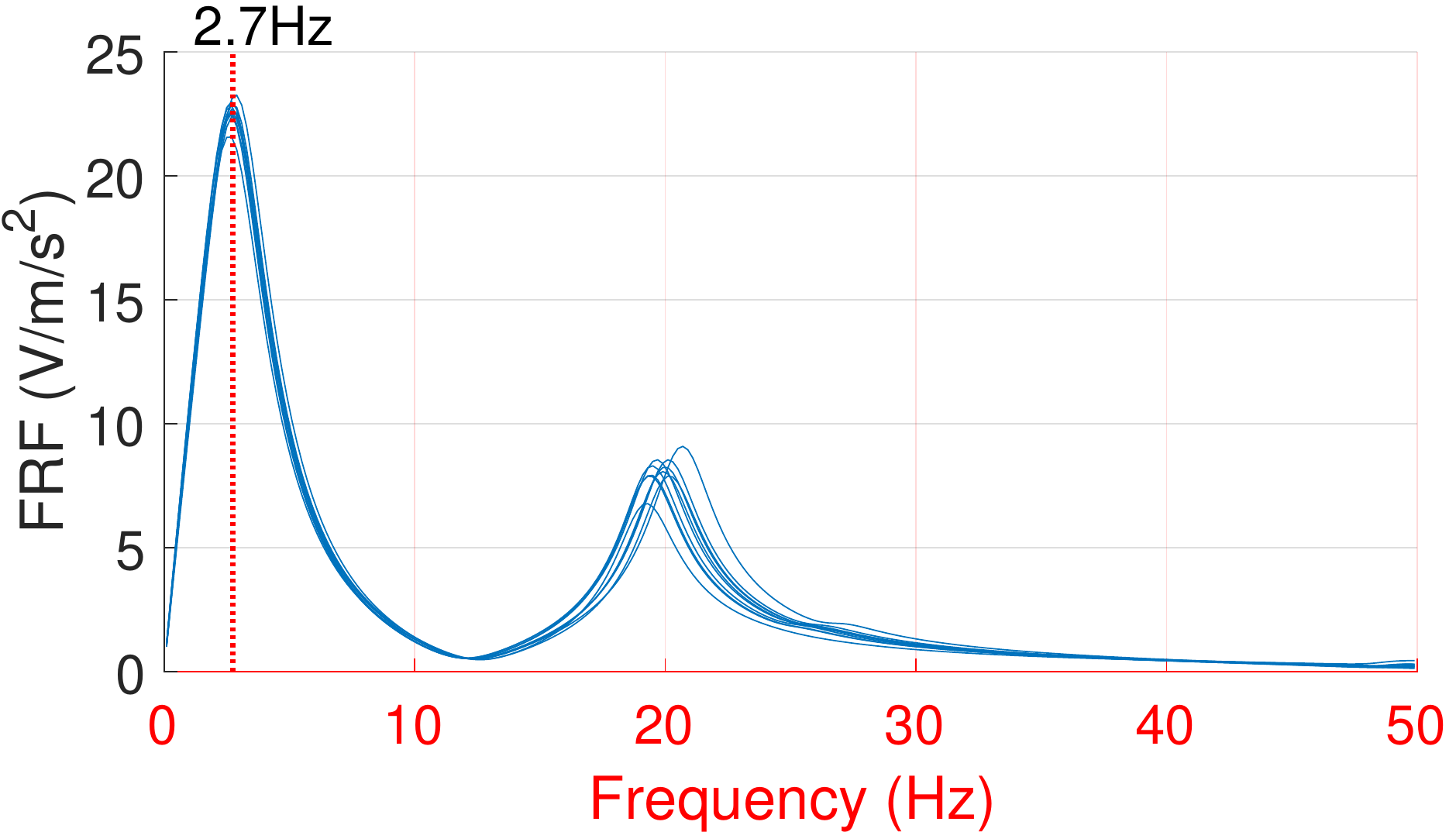}
       \caption{FRFs of all PEHs in cluster 2.}
    \end{subfigure}  
    \quad
    \begin{subfigure}[b]{0.28\textwidth}
        \includegraphics[width=\textwidth]{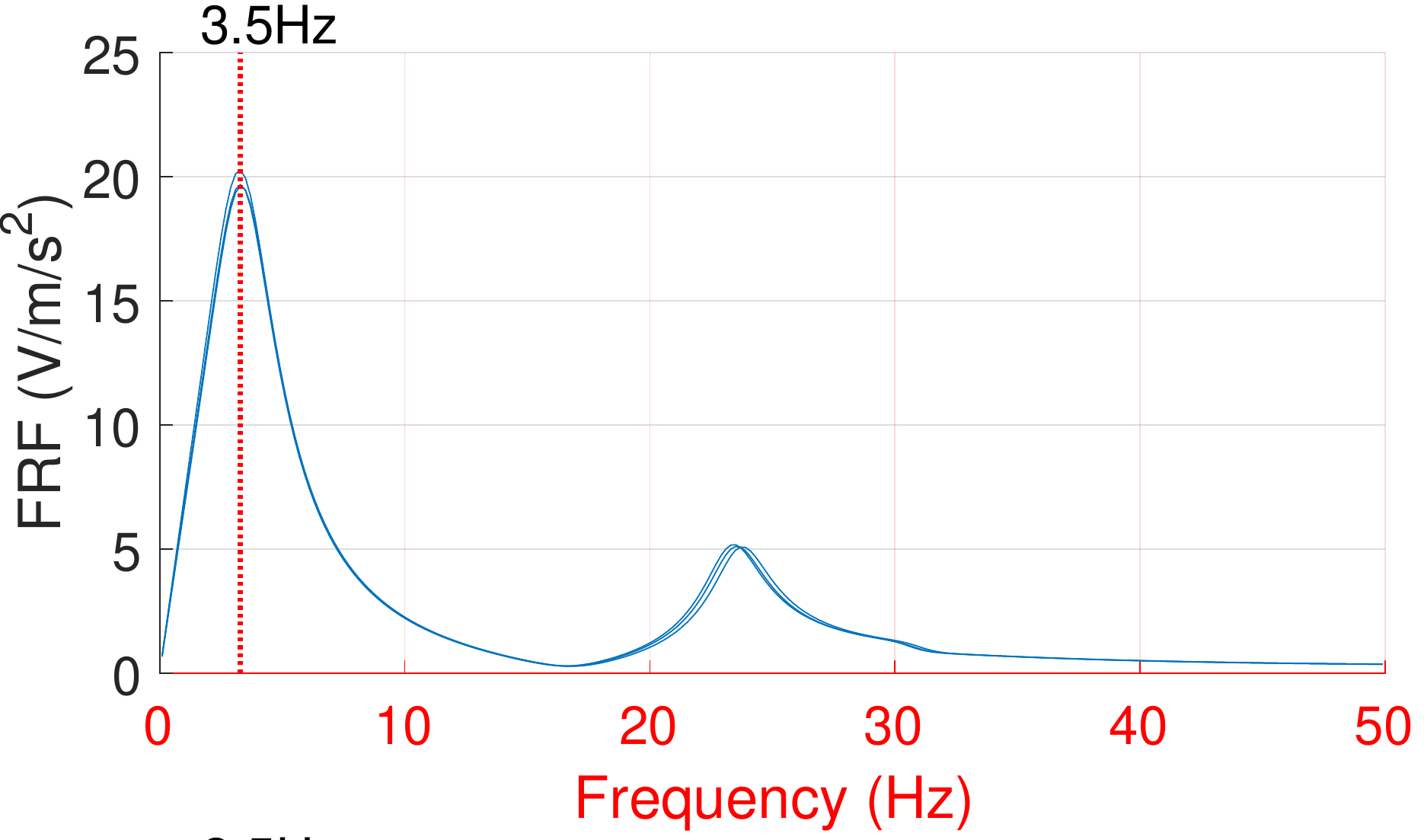}
        \caption{FRFs of all PEHs in cluster 3.}
    \end{subfigure}

    \begin{subfigure}[b]{0.28\textwidth}
    \vspace{0.5cm}
        \includegraphics[width=\textwidth]{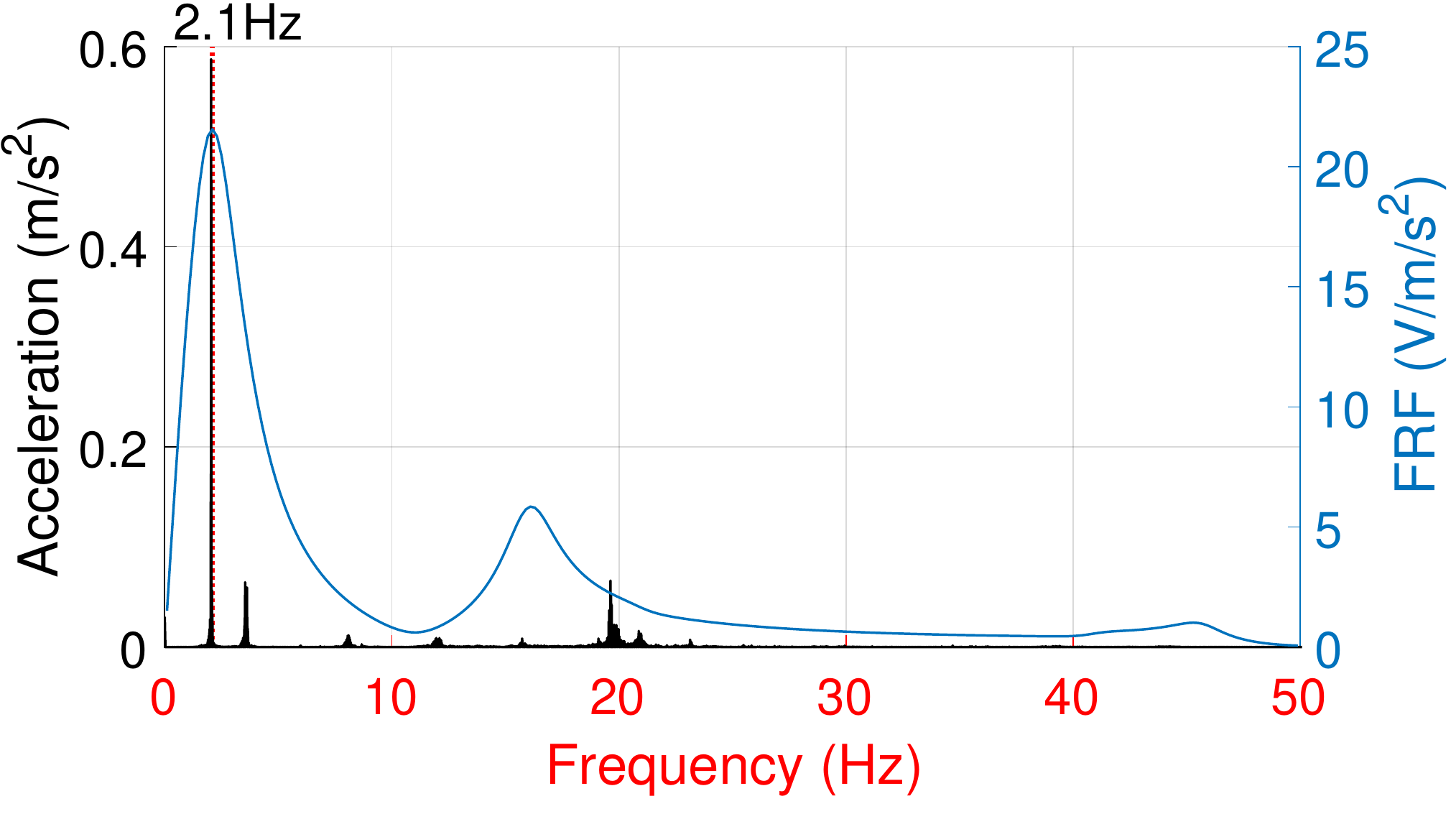}
       \caption{FRF of the centroid of cluster 1 and a sample bridge acceleration spectrum.}
    \end{subfigure} 
    \quad
    \begin{subfigure}[b]{0.28\textwidth}
    \vspace{0.5cm}
        \includegraphics[width=\textwidth]{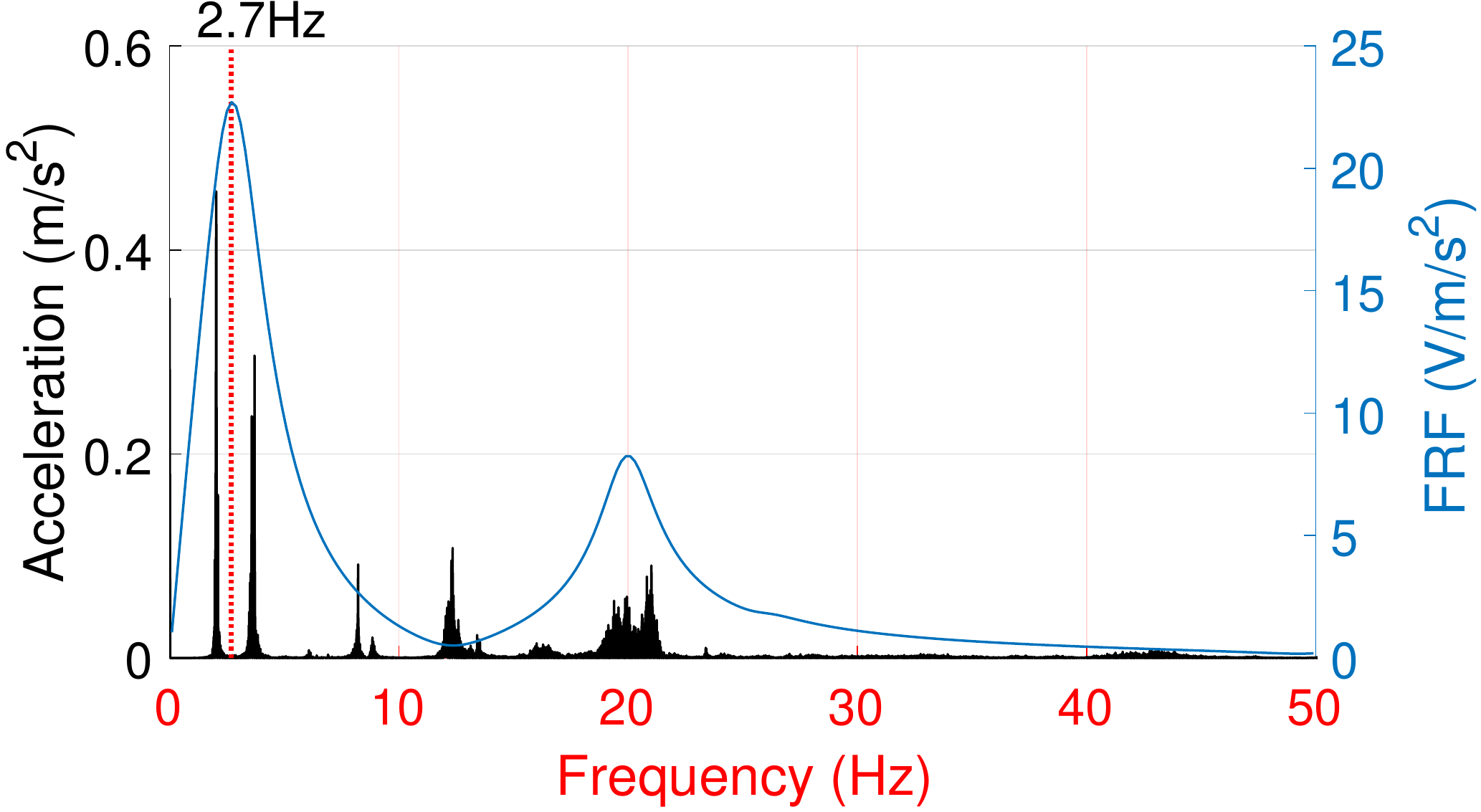}
        \caption{FRF of the centroid of cluster 2 and a sample bridge acceleration spectrum.}
    \end{subfigure}
    \quad
    \begin{subfigure}[b]{0.28\textwidth}
    \vspace{0.5cm}
        \includegraphics[width=\textwidth]{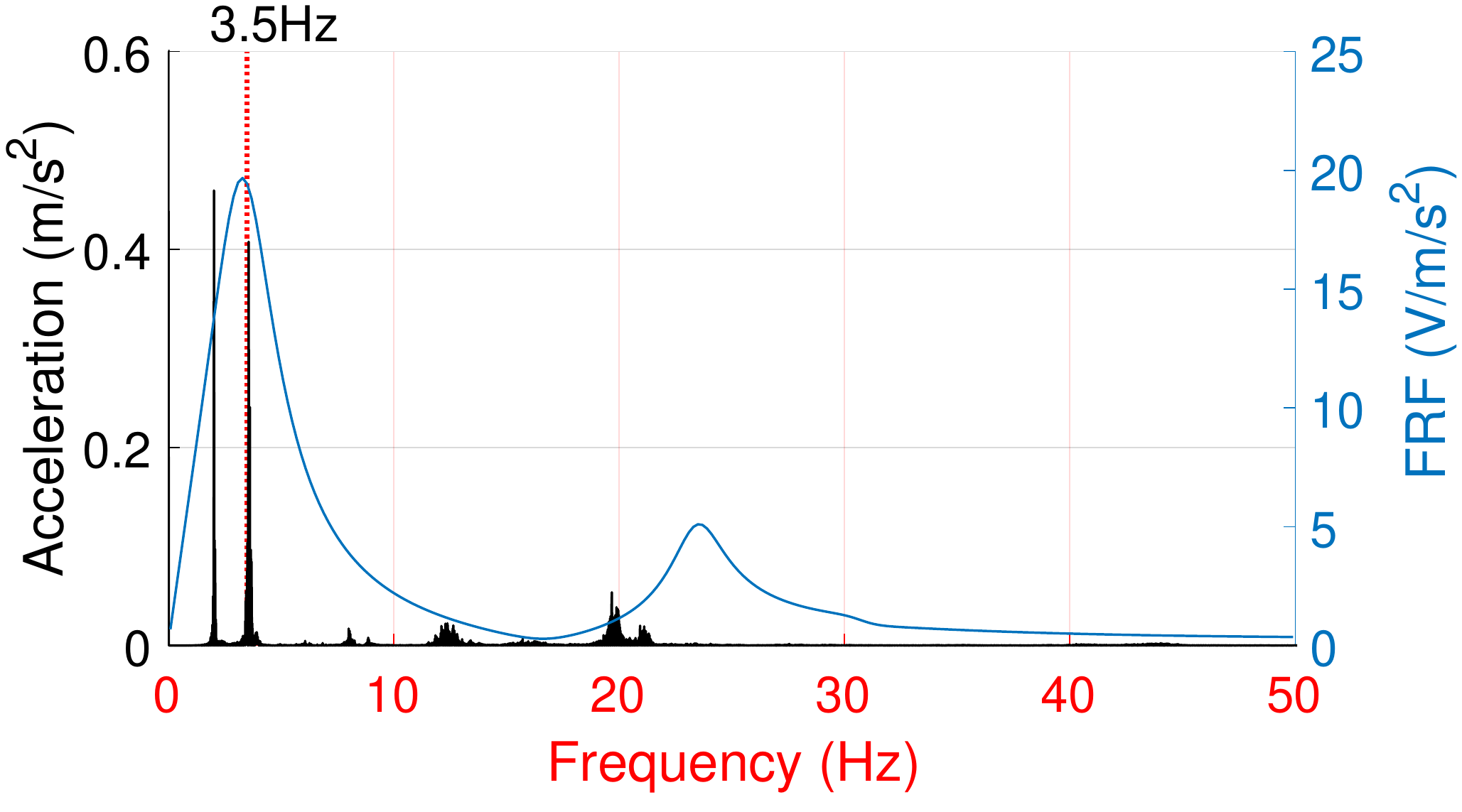}
       \caption{FRF of the centroid of cluster 3 and a sample bridge acceleration spectrum.}
    \end{subfigure}
 \caption{Illustration of the FRF response (blue line) of three candidates for sensor location A2 and their relationship with the Fourier acceleration spectrum of the input excitation (black line).}   
     \label{4.1.4}
\end{figure}

The Frequency Response Functions (FRFs) of the PEHs in clusters 1, 2, and 3 are shown in Figures \ref{4.1.4} (a, b, c) respectively. As it is expected, these figures indicate that the PEHs in the same cluster have similar frequency characteristics, with the natural frequencies being 2.1 Hz for cluster 1, 2.7 Hz for cluster 2, and 3.5 Hz for cluster 3. Additionally, Figures \ref{4.1.4} (d, e, f) display the FRFs of the cluster centroids and a sample Fourier spectrum of the input excitation signal from the same cluster. These plots reveal that the optimal fundamental frequency of the FRF aligns with the first peak of acceleration for cluster 1, lies between the first and second peaks of acceleration for cluster 2, and aligns with the second peak of acceleration for cluster 3. This result challenges the conventional approach to PEH design, which involves tuning the device to the fundamental frequency of the base structure. Instead, it shows that the optimal design depends on the input excitation signals. This analysis suggests that the variation in $L$ plays a crucial role in tuning the fundamental frequency of the device to the optimal harvesting frequency, while the variation in $l$ adjusts the energy absorbed by the device from the first two vibration modes, leading to improved energy output in some cases.


One crucial question that arises is the reason behind the different excitation signals observed during different hours of the day. These differences can be seen in the frequency spectrum of the input signals. In the next section, we will explore if there is a correlation between the volume of traffic and the dynamic response of the bridge, as amount of traffic is one of the key factors that change during the day.

\subsection{Effect of Traffic on Optimal Designs} \label{4.2.4}
In order to understand the PEH maximising the output efficiency, it is necessary to explore why the candidates generated under different inputs are similar. The input excitation window serves as an important variable affecting the PEH design during the optimisation process, so analyzing the excitation characteristics behind these PEH designs cannot be ignored. Since the traffic characteristics of vehicles are an important excitation source during the operation of bridges, this section takes the dynamic response of PEH at sensor location A2 as an example to explore how traffic affects energy generation and shape optimisation of devices by identifying the number of vehicles in consecutive windows.

\begin{figure}[H]
	\centering
	\includegraphics[width=0.95\linewidth]{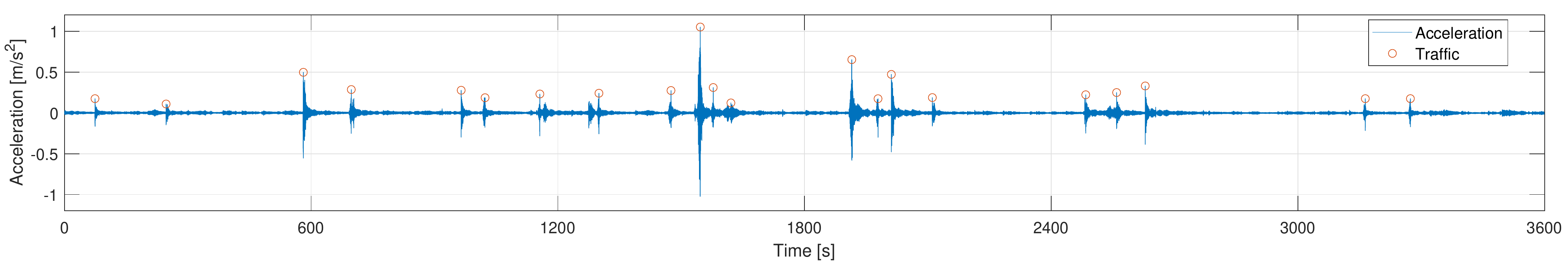}
	\caption{Number of vehicles identified in one-hour time window.}
	\label{4.3.2}
\end{figure}

Figure \ref{4.3.2} shows the number of vehicles identified from a typical acceleration time signal with a threshold of $0.2m/s^2$ of acceleration due to a passing vehicle. In order to analyse the impact of traffic volume on energy harvesting, the traffic in this case is divided into three levels according to the number of vehicles, as shown in Table \ref{T4.4}. As seen, the bridge is subject to a low traffic flow due to the nature of this bridge connecting two university campuses together and is mainly used by university staff/students. 

\begin{figure}[H]
	\centering
	\includegraphics[width=1\linewidth]{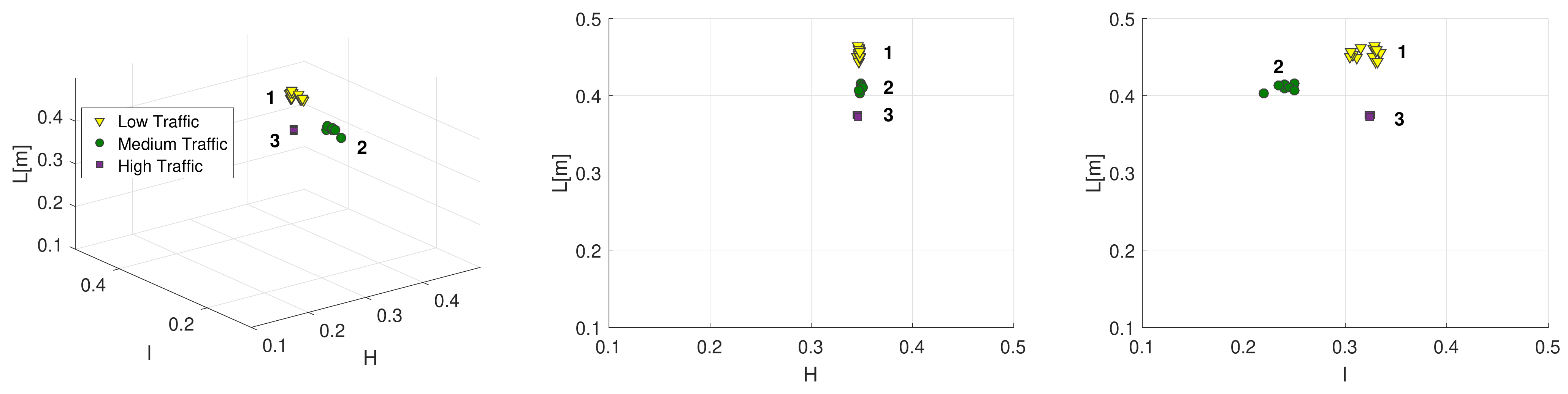}
	\caption{Three clusters of optimal PEHs corresponding to three traffic classes.}
	\label{Traffic}
\end{figure}

\begin{figure}[H]
	\centering
	\includegraphics[width=0.95\linewidth]{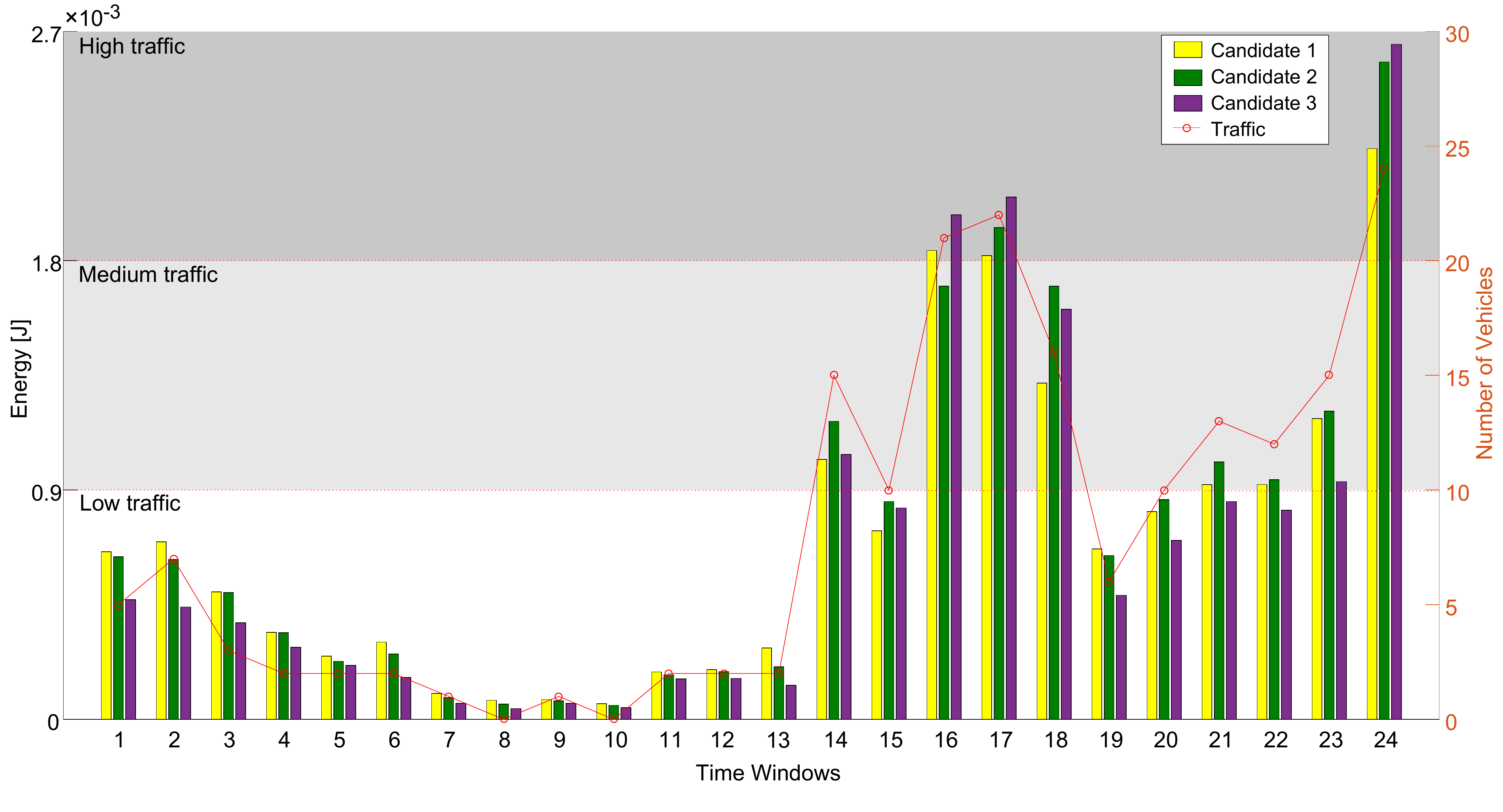}
	\caption{Energy harvested by three best candidates at sensor location A2 and its relationship to the number of vehicles.}
	\label{4.3.1}
\end{figure}

By classifying all the optimal PEHs at sensor location A2 according to the traffic intensity shown in Table \ref{T4.4}, it is found that the clustering results correspond to different traffic volumes (Figure \ref{Traffic}), explaining why the same location has different optimal geometric designs of PEH. Figure \ref{4.3.1} shows the dynamic responses of the optimal PEHs at sensor location A2 under different traffic excitation intensities, illustrating that the energy output is proportional to the traffic volume. We can see that for all windows with low traffic (1 - 13, 19) candidate 1 (i.e. design with the fundamental frequency 2.1 Hz) produces the most energy; for windows with medium traffic (14, 15, 18, 20 - 23), candidate 2 (with fundamental frequency 2.7 Hz) produces the most energy; for windows with high traffic (16, 17, 24) candidate 3 (with fundamental frequency 3.5 Hz) is the most efficient one. In other words, different types of best devices are matched to different traffic volumes. When there are very few passing vehicles, optimal PEH at sensor location A2 has the fundamental frequency that coincides with the first peak of input acceleration to generate the most energy; when the number of vehicles starts to increase, the fundamental frequency of the PEH starts to be affected by the second frequency; when the traffic volume is high, the maximum energy harvesting capability is achieved by tuning the device to the second frequency. 

\begin{table}[H]
\centering
\renewcommand\arraystretch{1.5}
\caption{Traffic classes of the study case}
\begin{tabular}{ccl}
\cline{1-2}
Traffic Class  & Number of Vehicles/hr &  \\ \cline{1-2}
Low Traffic    & 0-10                  &  \\
Medium Traffic & 10-20                 &  \\
High Traffic   & Above 20              &  \\ \cline{1-2}
\end{tabular}%
\label{T4.4}
\end{table}

\subsection{Candidates at all Locations} \label{4.2.2}
\begin{figure}[H]
	\centering
	\includegraphics[width=1\linewidth]{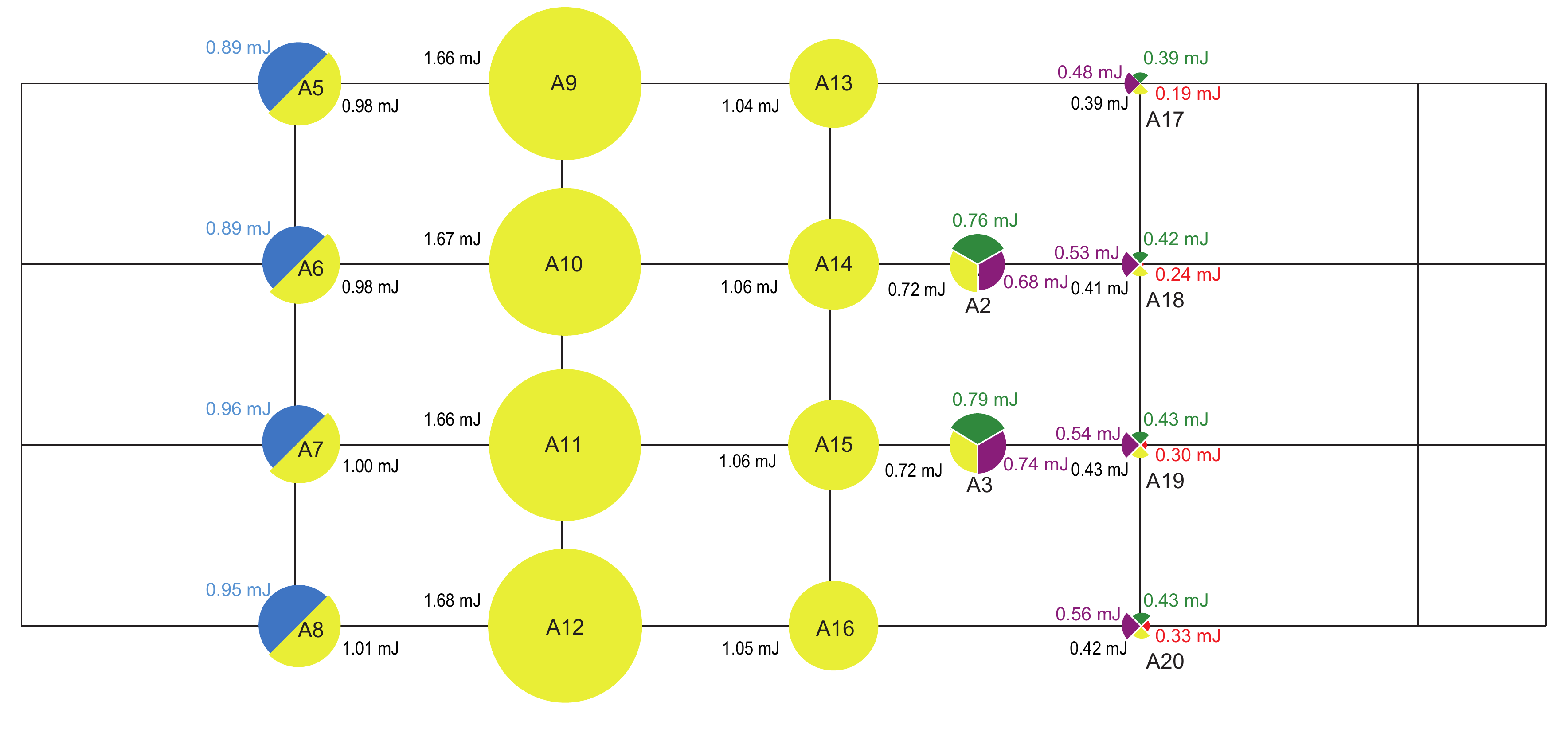}
	\caption{Energy harvesting capabilities of the candidates under a 24-hour continuous window at all locations.}
	\label{F4.2.2}
\end{figure}

The purpose of this research is to assess how the installation location impacts the performance of a piezoelectric energy harvester (PEH). The optimisation framework was applied to 17 additional positions, as described in section \ref{4.2.1}, to determine the best PEH design for each site. The results are presented in Figure \ref{F4.2.2}, where similar geometries are indicated by the same color-coding, and the size radius of each circular segment reflects the energy harvested from a 24-hour signal.

Five types of PEHs have been identified, each represented by a different color: yellow (candidate 1), green (candidate 2), purple (candidate 3), blue (candidate 4), and red (candidate 5). The results show that PEHs located in the same Cross-Girder (CG) direction produce a similar amount of energy (CG2: A5 - A8, CG3: A9 - A12, CG4: A13 - A16, CG5: A17 - A20 and A2 - A3). Additionally, the best candidate(s) for each location can be observed as follows: A9 - A16 have one optimal PEH, A5 - A8 have two optimal PEHs, A2 - A3 have three optimal PEHs, and A17 - A20 have four optimal PEHs.

\begin{figure}[H]
	\centering
	\includegraphics[width=0.95\linewidth]{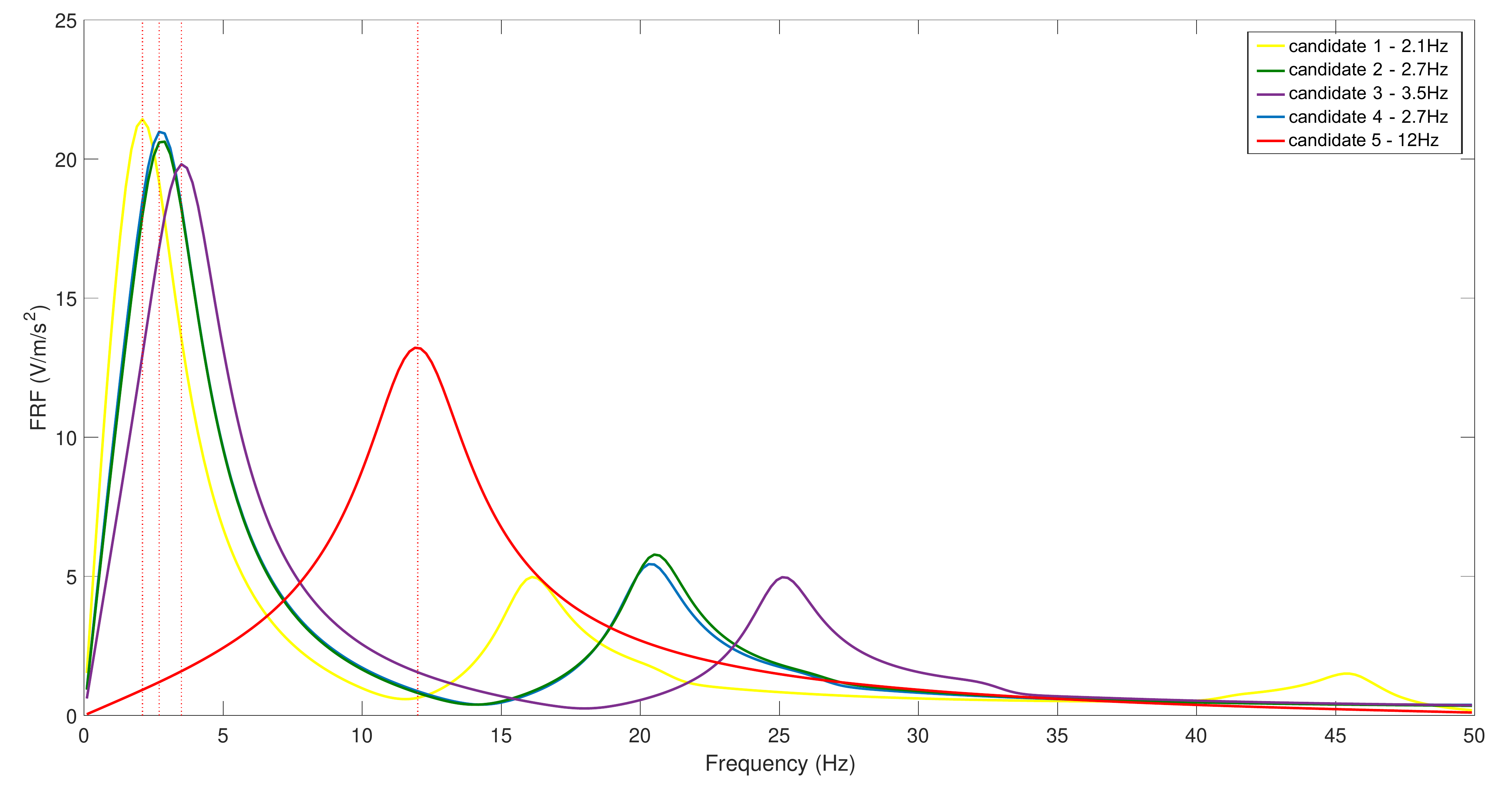}
	\caption{Illustration of the FRF features of the five types of optimal PEHs.}
	\label{A4.1.1}
\end{figure}

The FRFs of the best piezoelectric devices are displayed in Figure \ref{A4.1.1}. These devices have fundamental frequencies of 2.1 Hz (candidate 1), 2.7 Hz (candidates 2 and 4), 3.5 Hz (candidate 3), and 12 Hz (candidate 5). The variation in optimal PEH characteristics across different locations can be explained by the bridge bending modes, shown in Figure \ref{F4.2.3}. It can be seen that the geometric parameters of candidate 1 (yellow) under all time windows at the middle girder location (A9-A16) are centralised with a natural frequency close to the first vibration mode (2 Hz). This is because the middle girder location experiences large displacement in the first bending mode, resulting in the PEH at that position collecting more energy from that mode.

\begin{figure}[H]
	\centering
	\includegraphics[width=0.65\linewidth]{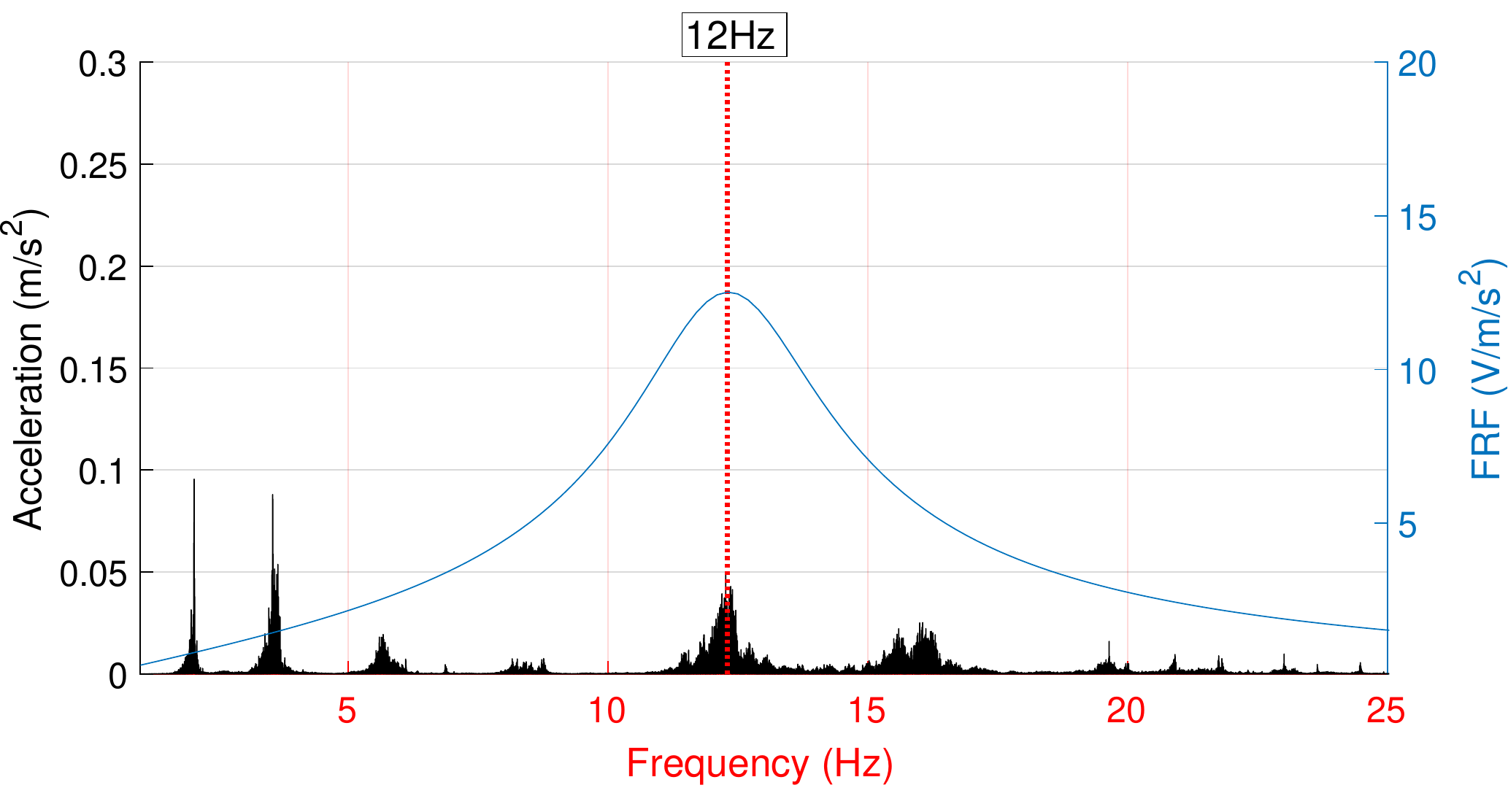}
	\caption{Illustration of the Fourier spectrum and the FRFs of candidates with a fundamental frequency of 12 Hz.}
	\label{12Hz}
\end{figure}

For other locations, however, the effects of other vibration modes are significant. As mentioned earlier, candidate 2 (green) and candidate 3 (purple) at A2 demonstrate the effect of the second deformation mode of the bridge on the PEH geometric design. From the conclusion in Section \ref{4.2.4}, this is because the second vibration mode is excited as the traffic volume increases, and A2 and A3 can harvest more energy from the second mode due to their peak positions. Therefore the optimal PEHs at A2 and A3 are optimised by tuning to the second frequency (candidate 3) when there are enough vehicles to fully excite the second mode. Note that increasing the amount of energy absorbed from one mode necessarily reduces the amount of energy absorbed from the other mode. This can be seen from the comparison of candidate 4 at A5-A8 (blue) and candidate 2 at A2-A3 (green) in Figure \ref{A4.1.1}, although they have the same fundamental frequency, they have different $l$. This is because the optimal devices at A2 and A3 increase the energy absorbed from the second mode by decreasing $l$. But it is not desirable to lose the energy available in the first mode at A5-A8 . Furthermore, in terms of A17-A20, the diversity of candidates is increased. Candidate 5 (red) near 12 Hz reveals the possibility of the PEH fundamental frequency tuned to the natural frequencies of other modes such as the widest acceleration range. This occurs when the first two acceleration peaks of the traffic excitation are not significant, as shown in Figure \ref{12Hz}.



Although it is concluded that the PEH configuration at different positions is influenced by the structural deformation mode, there are multiple candidates at each position due to the difference of input excitation signals. Therefore, it is necessary to identify the best location for energy harvesting and the corresponding best PEH design.

\subsection{Best Design and Location of PEH} \label{4.2.3}
Since the existence of multiple candidates at each location is not the best PEH design option, the comparison between the energy harvesting capabilities of the candidates under the 24-hour continuous acceleration signal input is performed at each location to summarize the design guidelines for the best PEHs that can maximise energy harvesting in a real environment based on bridge structural characteristics and traffic excitations. The result is shown in Figure \ref{F4.2.2}, where the color represents the PEH configuration, and the size of the circular segment is proportional to the amount of collected energy. 

The output responses presented in Figure \ref{F4.2.2} reveal the importance of location for energy harvesting. PEH designs in the optimal locations have the highest energy harvesting efficiency. For example, location A12 collects three times more energy than location A20 with the same equipment. So identifying the optimal locations is an important step in maximising power generation from the structure. The effect of location on energy harvesting is explained by the bending modes of the bridge deck (Figure \ref{F4.2.3}). It can be seen that the optimal locations for maximum power generation (A9, A10, A11, and A12) corresponds to where the maximum deflection occurs in the first deformation mode (near 2 Hz). Thus the PEH candidate 1 (yellow) has the potential to absorb the most energy from the first vibration mode. Therefore, the optimal installation locations for PEHs on bridges is revealed based on vibration modes that represent the dynamic behavior of the structure. It is worth noting that due to the low traffic volume in the operating environment of the bridge considered in this work, energy harvesting close to the first natural frequency is the optimal choice for this case. For example, although A17-A20 has a large peak in the second bending mode, it is the location with the worst energy harvesting ability due to the excitation characteristics of the input signal. 

\begin{figure}[H]
	\centering
	\includegraphics[width=1\linewidth]{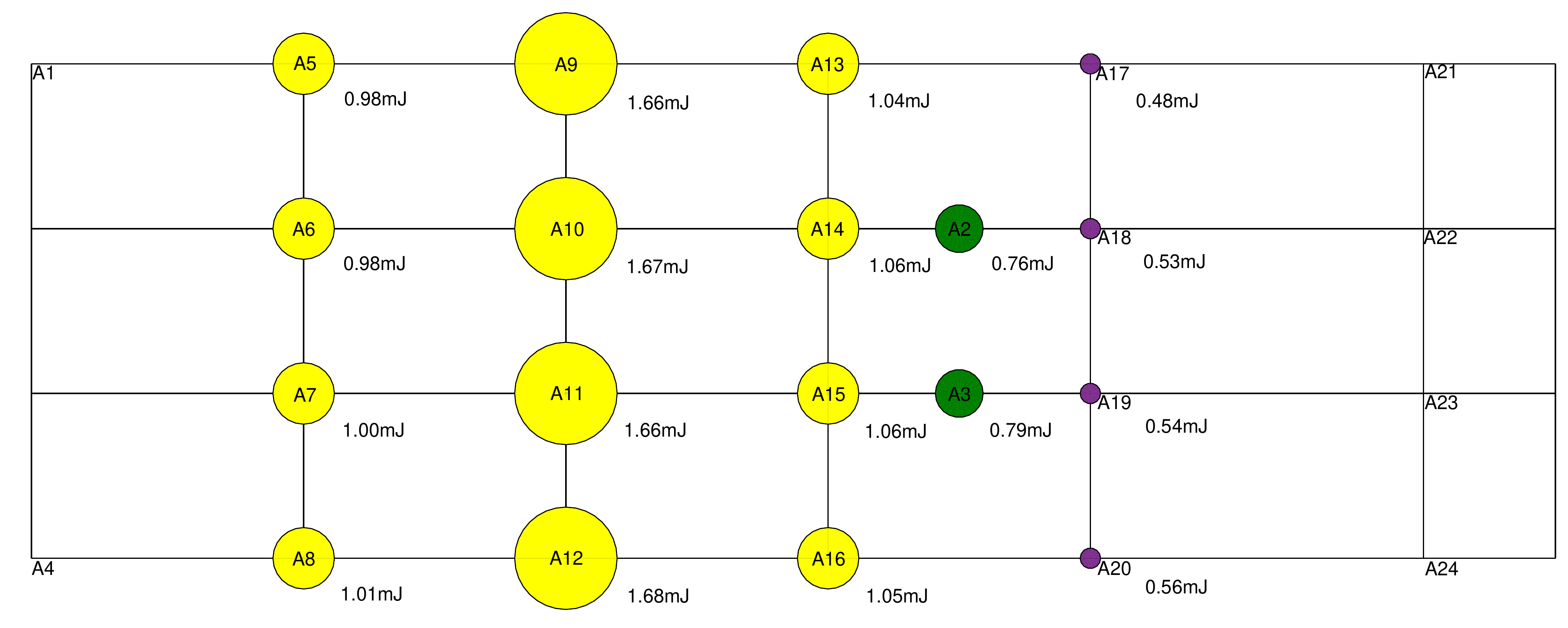}
	\caption{Energy harvesting capabilities of the three final candidates (denoted by yellow, green, and purple colors) under a 24-hour continuous excitation time window.}
	\label{F4.2.5}
\end{figure}
\begin{figure}[H]
	\centering
	\includegraphics[width=1\linewidth]{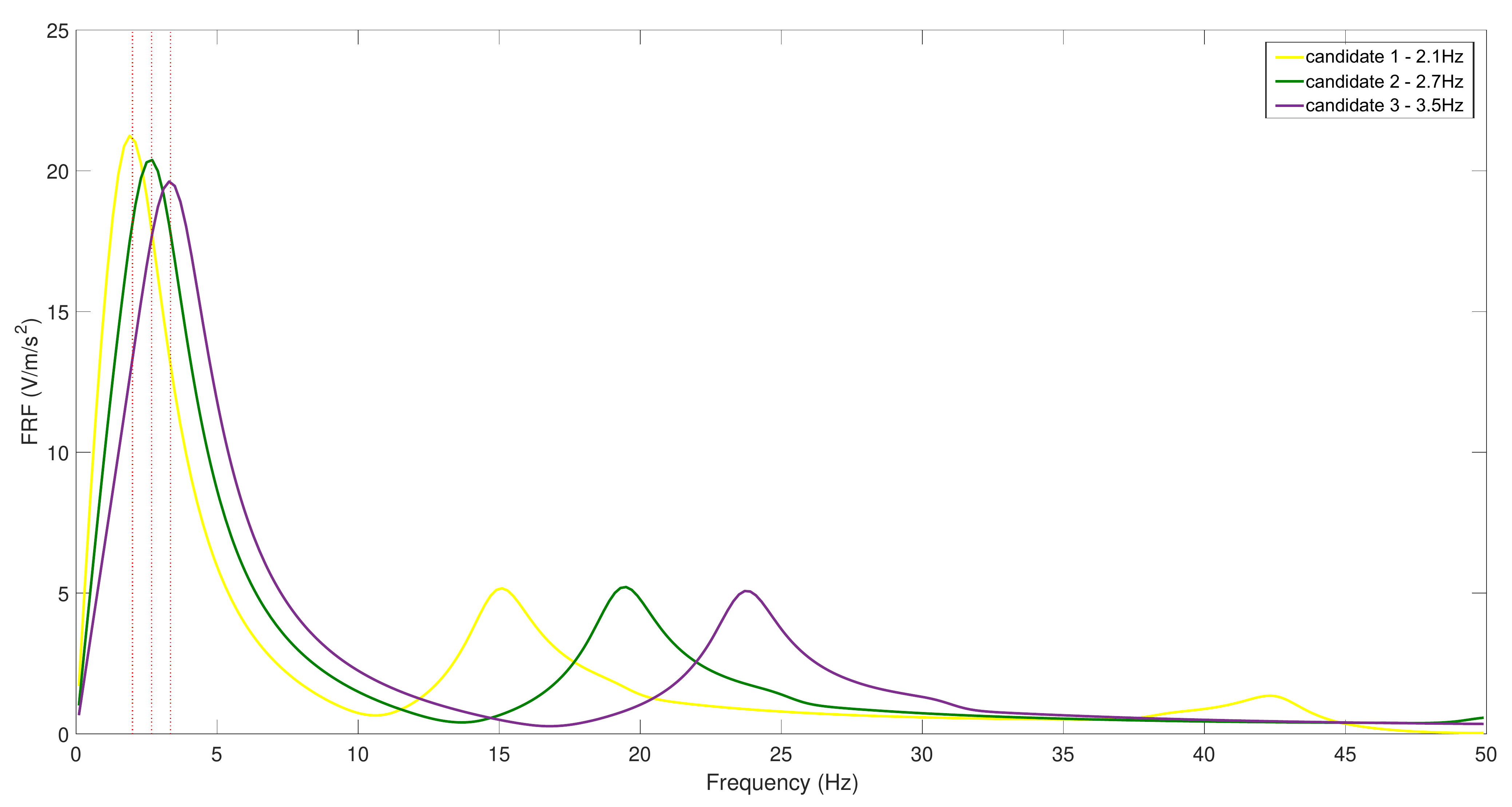}
	\caption{Illustration of the FRF response of the three final devices.}
	\label{F4.2.6}
\end{figure}

In addition, the best options with the largest output capability at each position under the 24-hour acceleration signal are found, and the types of best PEHs on the entire structure are reduced to 3, as shown in Figure \ref{F4.2.5}. Figure \ref{F4.2.6} shows the FRFs and the corresponding fundamental frequencies for the three final devices, compared with the mode shapes (Figure \ref{F4.2.3}) to conclude that the best design at the position with large deformation in the first bending mode tends to have similar frequency characteristics to the first vibration mode. Therefore, a consistent phenomenon is exhibited for the optimised designs from A5 to A16. For the optimal design at the edge girders (A17, A18, A19 and A20), little energy is harvested from the first vibration mode due to the small displacement of the first bending mode shown in Figure \ref{F4.2.3}b. Thus the PEHs mounted at the edges can absorb most of the energy from the second mode through geometry optimisation, i.e. tuning to the second natural frequency.

The best candidates at A2 and A3 indicate that the fundamental frequency of the PEH device can not only improve energy harvesting efficiency by matching the amplitude peaks of the acceleration spectrum but also can be tuned to a certain frequency range to harvest energy from multiple vibration modes, which depends on the deformation mode of the structure and the strength of the external excitation. Although A2 and A3 have high displacements in the second mode, due to the characteristics of the traffic flow under the real operating conditions of the cable-stayed bridge in this work, the first mode is more inclined to be excited, and the first acceleration peak is higher than the second peak. Therefore, the energy harvested from the first mode is non-negligible, the PEH is designed to harvest energy from the first two modes by synchronising the natural frequencies between the first two frequencies.


\begin{figure}[H]
	\centering
	\includegraphics[width=1\linewidth]{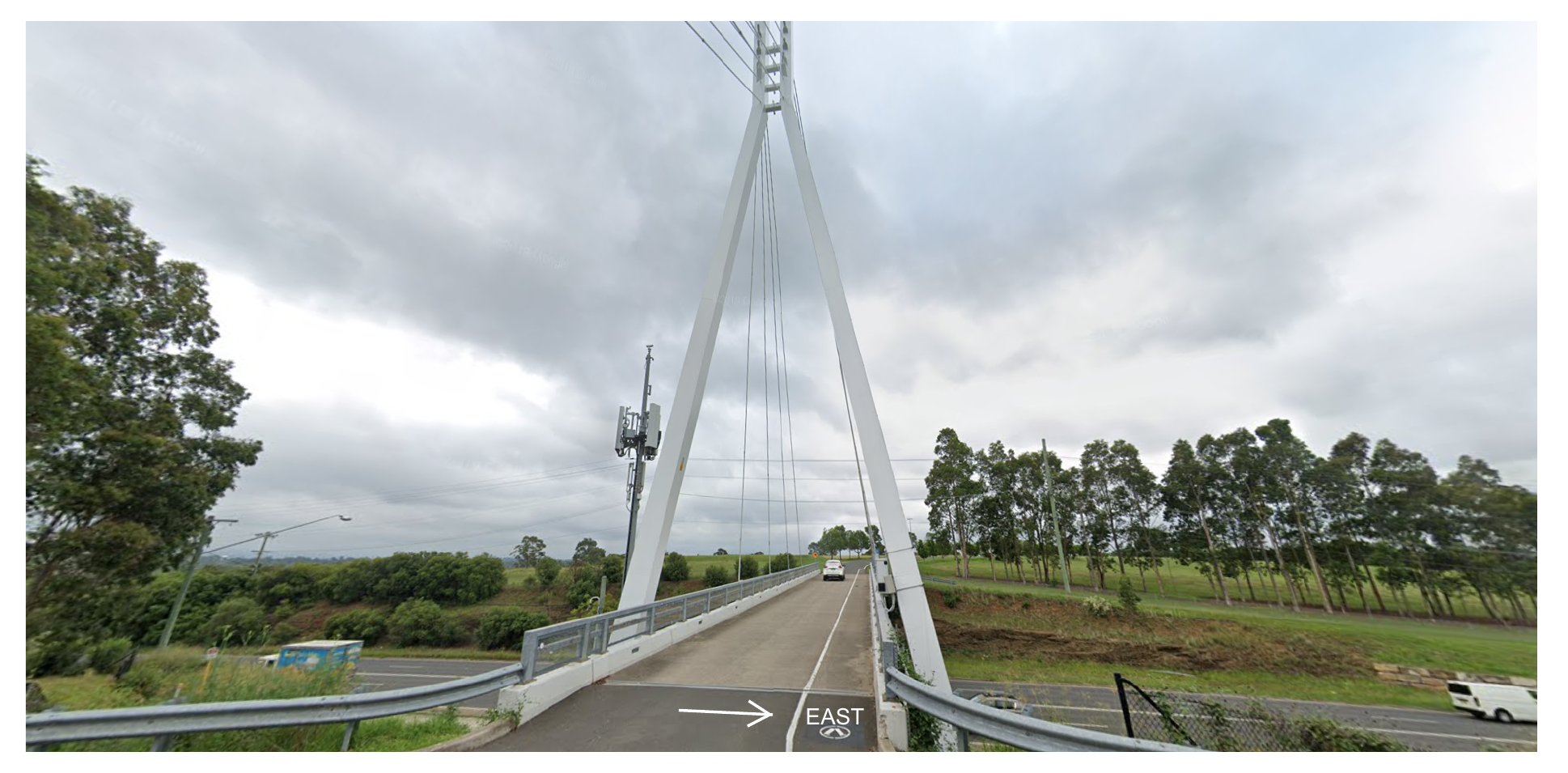}
	\caption{Bridge Lane Site Map.}
	\label{F4.2.4}
\end{figure}

Another interesting finding is that the harvesters placed on the west side of this structure always have a higher energy output than the east side (Figure \ref{F4.2.2}), which is explained by the traffic lane shown in the bridge site map (Figure \ref{F4.2.4}). It is clear that the driveway (west) has a stronger excitation than the pedestrian lane (east), so locations with stronger excitation collect more energy.

Therefore, the above optimal locations for energy harvesting and corresponding optimal designs demonstrate that the framework can be applied to identify the geometric configuration and installation position of a limited number of PEHs on the bridge to achieve the best energy harvesting performance. At the same time, this paper also summarises the design rules of piezoelectric devices to maximise the energy harvesting on the structure by considering the influence of location and traffic.

%% file: Sections/conclusion.tex
This work comprehensively investigated how the installation location of the piezoelectric energy harvester affects the energy output efficiency based on the vibration excitation input in the real environment provided by a cable-stayed bridge, and proposed a PEH shape optimisation framework that maximises energy harvesting with structural location as a variable. In pursuit of the PEH design and installation location that maximises harvested energy, we extended our previous PEH optimisation framework to consider time windows more consistent with vibrational excitations in natural environments. The results showed that the transformable energy at different locations varies greatly, especially along the length of the bridge (traffic direction), which has a difference from 0.42 mJ (A20) to 1.68 mJ (A12), an increase of 3 times.

With the implementation of the proposed optimisation framework on the entire bridge, the optimal placement of the PEH was determined. To highlight the impact of location on energy harvesting, the harvesting capacity of the best PEH at various positions was correlated with the deformation modes of the bridge deck. The results showed that the location with the maximum deformation under the first vibration mode was the most suitable for energy harvesting, and the best PEH had a fundamental frequency close to the first frequency.

We also explored the impact of different levels of traffic intensity on energy harvesting by grouping traffic volumes. The results showed differences in available energy and optimal PEH shape under different traffic conditions. Firstly, the convertible mechanical energy varies greatly under different traffic volumes. In addition, the optimal PEH designs identified under different traffic volumes are different. This is because the excitation degree of different traffic has a relation to the vibration mode of the bridge, and then guides the fundamental frequency of the device to move within a certain frequency range to seek the maximisation of energy output.

In summary, this work realised the practical optimisation of PEH performance in the spatial domain on a real structure. And it was proposed to design the optimal shape and installation location (predominantly its variation along the length of the bridge) of a PEH according to the dominant traffic volume and the deformation mode in the operation of the structure to obtain the maximum energy harvesting capacity.

An interesting finding is that the energy harvested at the traffic lane of the bridge is higher than that at the sidewalk, which inspired us to further explore the impact of more detailed traffic characteristics on energy harvesting, not just the number of vehicles, but also the driving speed and direction, etc.

%% file: Sections/Appendix.tex
\begin{equation}
    M_{IJ} = \int_{\Omega _s} \rho_s(N_IN_J+z^2N_{I,x}N_{J,x}+z^2N_{I,y}N_{J,y}) \text{d}\Omega _s + 2\int_{\Omega _p} \rho_p(N_IN_J+z^2N_{I,x}N_{J,x}+z^2N_{I,y}N_{J,y}) \text{d}\Omega _p 
    \label{A1}
\end{equation}

\begin{equation}
    K_{IJ} = \int_{\Omega _s} z^2\mathbf{B}_{I}^{T}\mathbf{c}_s\mathbf{B}_J \text{d}\Omega _s + 2\int_{\Omega _s} z^2\mathbf{B}_{I}^{T}\mathbf{c}_p^E\mathbf{B}_J \text{d}\Omega _s 
    \label{e8}
\end{equation}

\begin{equation}
    \Theta _{I} = \int_{\Omega _p} z^2\mathbf{B}_{I}^{T}\mathbf{e}^T\mathbf{Z} \text{d}\Omega _p
    \label{e9}
\end{equation}

\begin{equation}
\begin{aligned}
    F_{I} = \int_{\Omega _s} \rho_sN_I \text{d}\Omega _s+\int_{\Omega _p} \rho_pN_I \text{d}\Omega _p
+\sum_{J=1}^{k}\int_{\Omega _s} \rho_sz^2(N_{I,x}N_{J,x}+N_{I,y}N_{J,y}) \text{d}\Omega _s \\+ \sum_{J=1}^{k}\int_{\Omega _p} \rho_pz^2(N_{I,x}N_{J,x}+N_{I,y}N_{J,y}) \text{d}\Omega _p
    \label{pluse9}
    \end{aligned}
\end{equation}

\begin{equation}
   \mathbf{B}_I = \left \{ -N_{I,xx} \quad -N_{I,yy} \quad -2N_{I,xy} \right \}^T
    \label{e10}
\end{equation}

\begin{equation}
   \mathbf{Z} = \left \{ 0 \quad 0 \quad \frac{1}{h_p} \right \}^T
    \label{e11}
\end{equation}
where $\rho_p$ and $\rho_s$ are the density of the piezoelectric material and the base, respectively.

%% file: main.bbl
\begin{thebibliography}{30}
\expandafter\ifx\csname natexlab\endcsname\relax\def\natexlab#1{#1}\fi
\providecommand{\bibinfo}[2]{#2}
\ifx\xfnm\relax \def\xfnm[#1]{\unskip,\space#1}\fi
\bibitem[{Brownjohn(2006)}]{I1}
\bibinfo{author}{J.~M.~W. Brownjohn},
\newblock \bibinfo{title}{Structural health monitoring of civil
  infrastructure},
\newblock \bibinfo{journal}{Philosophical Transactions of the Royal Society A:
  Mathematical, Physical and Engineering Sciences} \bibinfo{volume}{365}
  (\bibinfo{year}{2006}) \bibinfo{pages}{589 -- 622}.
\bibitem[{Heidemann et~al.(2012)Heidemann, Stojanovic, and Zorzi}]{27}
\bibinfo{author}{J.~Heidemann}, \bibinfo{author}{M.~Stojanovic},
  \bibinfo{author}{M.~Zorzi},
\newblock \bibinfo{title}{Underwater sensor networks: applications, advances
  and challenges.},
\newblock \bibinfo{journal}{Philosophical Transactions} \bibinfo{volume}{370}
  (\bibinfo{year}{2012}) \bibinfo{pages}{158}.
\bibitem[{Asadollahi and Jian(2017)}]{I2}
\bibinfo{author}{P.~Asadollahi}, \bibinfo{author}{L.~Jian},
\newblock \bibinfo{title}{Statistical analysis of modal properties of a
  cable-stayed bridge through long-term structural health monitoring with
  wireless smart sensor networks},
\newblock \bibinfo{journal}{Journal of Bridge Engineering} \bibinfo{volume}{22}
  (\bibinfo{year}{2017}) \bibinfo{pages}{04017051}.
\bibitem[{Elvin et~al.(2003)Elvin, Elvin, and Choi}]{70}
\bibinfo{author}{N.~Elvin}, \bibinfo{author}{A.~Elvin}, \bibinfo{author}{D.~H.
  Choi},
\newblock \bibinfo{title}{{A self-powered damage detection sensor}},
\newblock \bibinfo{journal}{Strain Analysis for Engineering Design}
  \bibinfo{volume}{38} (\bibinfo{year}{2003}) \bibinfo{pages}{115--124}.
\bibitem[{Lee et~al.(2003)Lee, Hsu, Hsiao, and Wu}]{54}
\bibinfo{author}{C.~K. Lee}, \bibinfo{author}{Y.~H. Hsu},
  \bibinfo{author}{W.~H. Hsiao}, \bibinfo{author}{J.~Wu},
\newblock \bibinfo{title}{Electrical and mechanical field interactions of
  piezoelectric systems: foundation of smart structures-based piezoelectric
  sensors and actuators, and free-fall sensors},
\newblock \bibinfo{journal}{Smart Materials \& Structures} \bibinfo{volume}{13}
  (\bibinfo{year}{2003}) \bibinfo{pages}{1090}.
\bibitem[{Khalid et~al.(2015)Khalid, Redhewal, Kumar, and Srivastav}]{I3}
\bibinfo{author}{A.~Khalid}, \bibinfo{author}{A.~K. Redhewal},
  \bibinfo{author}{M.~Kumar}, \bibinfo{author}{A.~Srivastav},
\newblock \bibinfo{title}{Piezoelectric vibration harvesters based on
  vibrations of cantilevered bimorphs: A review},
\newblock \bibinfo{journal}{Materials Sciences \& Applications}
  \bibinfo{volume}{6} (\bibinfo{year}{2015}) \bibinfo{pages}{818--827}.
\bibitem[{Benasciutti et~al.(2010)Benasciutti, Moro, Zelenika, and Brusa}]{126}
\bibinfo{author}{D.~Benasciutti}, \bibinfo{author}{L.~Moro},
  \bibinfo{author}{S.~Zelenika}, \bibinfo{author}{E.~Brusa},
\newblock \bibinfo{title}{{Vibration energy scavenging via piezoelectric
  bimorphs of optimized shapes}} \bibinfo{volume}{16} (\bibinfo{year}{2010})
  \bibinfo{pages}{657--668}.
\bibitem[{Roundy et~al.(2005)Roundy, Leland, Baker, Carleton, Reilly, Lai,
  Otis, Rabaey, Wright, and Sundararajan}]{125}
\bibinfo{author}{S.~Roundy}, \bibinfo{author}{E.~S. Leland},
  \bibinfo{author}{J.~Baker}, \bibinfo{author}{E.~Carleton},
  \bibinfo{author}{E.~Reilly}, \bibinfo{author}{E.~Lai},
  \bibinfo{author}{B.~Otis}, \bibinfo{author}{J.~M. Rabaey},
  \bibinfo{author}{P.~K. Wright}, \bibinfo{author}{V.~Sundararajan},
\newblock \bibinfo{title}{Improving power output for vibration-based energy
  scavengers},
\newblock \bibinfo{journal}{IEEE Pervasive Computing} \bibinfo{volume}{4}
  (\bibinfo{year}{2005}) \bibinfo{pages}{28--36}.
\bibitem[{Peralta et~al.(2020)Peralta, Ruiz, Natarajan, and Atroshchenko}]{97}
\bibinfo{author}{P.~Peralta}, \bibinfo{author}{R.~Ruiz},
  \bibinfo{author}{S.~Natarajan}, \bibinfo{author}{E.~Atroshchenko},
\newblock \bibinfo{title}{Parametric study and shape optimization of
  piezoelectric energy harvesters by isogeometric analysis and kriging
  metamodeling},
\newblock \bibinfo{journal}{Journal of Sound and Vibration}
  \bibinfo{volume}{484} (\bibinfo{year}{2020}).
\bibitem[{Plessis et~al.(2005)Plessis, Huigsloot, and Discenzo}]{71}
\bibinfo{author}{A.~Plessis}, \bibinfo{author}{M.~J. Huigsloot},
  \bibinfo{author}{F.~D. Discenzo},
\newblock \bibinfo{title}{Resonant packaged piezoelectric power harvester for
  machinery health monitoring},
\newblock \bibinfo{journal}{Proceedings of {SPIE} - The International Society
  for Optical Engineering} \bibinfo{volume}{5762} (\bibinfo{year}{2005})
  \bibinfo{pages}{224--235}.
\bibitem[{Farinholt et~al.(2010)Farinholt, Miller, Sifuentes, MacDonald, Park,
  and Farrar}]{a1}
\bibinfo{author}{K.~M. Farinholt}, \bibinfo{author}{N.~Miller},
  \bibinfo{author}{W.~Sifuentes}, \bibinfo{author}{J.~MacDonald},
  \bibinfo{author}{G.~Park}, \bibinfo{author}{C.~R. Farrar},
\newblock \bibinfo{title}{Energy harvesting and wireless energy transmission
  for embedded shm sensor nodes},
\newblock \bibinfo{journal}{Structural Health Monitoring} \bibinfo{volume}{9}
  (\bibinfo{year}{2010}) \bibinfo{pages}{269--280}.
\bibitem[{Cahill et~al.(2018)Cahill, Hazra, Karoumi, Mathewson, and
  Pakrashi}]{a2}
\bibinfo{author}{P.~Cahill}, \bibinfo{author}{B.~Hazra},
  \bibinfo{author}{R.~Karoumi}, \bibinfo{author}{A.~Mathewson},
  \bibinfo{author}{V.~Pakrashi},
\newblock \bibinfo{title}{Vibration energy harvesting based monitoring of an
  operational bridge undergoing forced vibration and train passage},
\newblock \bibinfo{journal}{Mechanical Systems and Signal Processing}
  \bibinfo{volume}{106} (\bibinfo{year}{2018}) \bibinfo{pages}{265--283}.
\bibitem[{Laura et~al.(1974)Laura, Pombo, and Susemihl}]{75}
\bibinfo{author}{P.~Laura}, \bibinfo{author}{J.~L. Pombo},
  \bibinfo{author}{E.~A. Susemihl},
\newblock \bibinfo{title}{A note on the vibrations of a clamped-free beam with
  a mass at the free end},
\newblock \bibinfo{journal}{Journal of Sound and Vibration}
  \bibinfo{volume}{37} (\bibinfo{year}{1974}) \bibinfo{pages}{161--168}.
\bibitem[{Erturk and Inman(2008)}]{77}
\bibinfo{author}{A.~Erturk}, \bibinfo{author}{D.~J. Inman},
\newblock \bibinfo{title}{A distributed parameter electromechanical model for
  cantilevered piezoelectric energy harvesters},
\newblock \bibinfo{journal}{Journal of Vibration and Acoustics}
  \bibinfo{volume}{130} (\bibinfo{year}{2008}) \bibinfo{pages}{1257--1261}.
\bibitem[{Crandall and Balise(1970)}]{79}
\bibinfo{author}{S.~H. Crandall}, \bibinfo{author}{P.~L. Balise},
\newblock \bibinfo{title}{Dynamics and mechanical and electromechanical
  systems},
\newblock \bibinfo{journal}{Phys. Today} \bibinfo{volume}{23}
  (\bibinfo{year}{1970}) \bibinfo{pages}{75--75}.
\bibitem[{Lu et~al.(2003)Lu, Lee, and Lim}]{81}
\bibinfo{author}{F.~Lu}, \bibinfo{author}{H.~P. Lee}, \bibinfo{author}{S.~P.
  Lim},
\newblock \bibinfo{title}{Modeling and analysis of micro piezoelectric power
  generators for micro-electromechanical-systems applications},
\newblock \bibinfo{journal}{Smart Materials and Structures}
  \bibinfo{volume}{13} (\bibinfo{year}{2003}) \bibinfo{pages}{57--63}.
\bibitem[{{Calderon Hurtado} et~al.(2022){Calderon Hurtado}, Peralta, Ruiz,
  {Makki Alamdari}, and Atroshchenko}]{extra1}
\bibinfo{author}{A.~{Calderon Hurtado}}, \bibinfo{author}{P.~Peralta},
  \bibinfo{author}{R.~Ruiz}, \bibinfo{author}{M.~{Makki Alamdari}},
  \bibinfo{author}{E.~Atroshchenko},
\newblock \bibinfo{title}{Shape optimization of piezoelectric energy harvesters
  of variable thickness},
\newblock \bibinfo{journal}{Journal of Sound and Vibration}
  \bibinfo{volume}{517} (\bibinfo{year}{2022}) \bibinfo{pages}{116503}.
\bibitem[{Peralta-Braz et~al.(2022)Peralta-Braz, Alamdari, Ruiz, Atroshchenko,
  and Hassan}]{ex35}
\bibinfo{author}{P.~Peralta-Braz}, \bibinfo{author}{M.~M. Alamdari},
  \bibinfo{author}{R.~O. Ruiz}, \bibinfo{author}{E.~Atroshchenko},
  \bibinfo{author}{M.~Hassan},
\newblock \bibinfo{title}{Design optimisation of piezoelectric energy
  harvesters for bridge infrastructure}  (\bibinfo{year}{2022}).
\bibitem[{Zhang et~al.(2018)Zhang, Xiang, Shi, and Zhan}]{ex40}
\bibinfo{author}{Z.~Zhang}, \bibinfo{author}{H.~Xiang},
  \bibinfo{author}{Z.~Shi}, \bibinfo{author}{J.~Zhan},
\newblock \bibinfo{title}{Experimental investigation on piezoelectric energy
  harvesting from vehicle-bridge coupling vibration},
\newblock \bibinfo{journal}{Energy Conversion and Management}
  \bibinfo{volume}{163} (\bibinfo{year}{2018}) \bibinfo{pages}{169--179}.
\bibitem[{Karimi et~al.(2016)Karimi, Karimi, Tikani, and Ziaei-Rad}]{4.3}
\bibinfo{author}{M.~Karimi}, \bibinfo{author}{A.~H. Karimi},
  \bibinfo{author}{R.~Tikani}, \bibinfo{author}{S.~Ziaei-Rad},
\newblock \bibinfo{title}{Experimental and theoretical investigations on
  piezoelectric-based energy harvesting from bridge vibrations under travelling
  vehicles},
\newblock \bibinfo{journal}{International Journal of Mechanical Sciences}
  (\bibinfo{year}{2016}) \bibinfo{pages}{1--11}.
\bibitem[{Song(2019)}]{ex41}
\bibinfo{author}{Y.~Song},
\newblock \bibinfo{title}{Finite-element implementation of piezoelectric energy
  harvesting system from vibrations of railway bridge},
\newblock \bibinfo{journal}{Journal of Energy Engineering}
  \bibinfo{volume}{145} (\bibinfo{year}{2019})
  \bibinfo{pages}{04018076.1--04018076.18}.
\bibitem[{Alamdari et~al.(2019)Alamdari, Kildashti, Samali, and Goudarzi}]{4.2}
\bibinfo{author}{M.~M. Alamdari}, \bibinfo{author}{K.~Kildashti},
  \bibinfo{author}{B.~Samali}, \bibinfo{author}{H.~V. Goudarzi},
\newblock \bibinfo{title}{Damage diagnosis in bridge structures using rotation
  influence line: Validation on a cable-stayed bridge},
\newblock \bibinfo{journal}{Engineering Structures} \bibinfo{volume}{185}
  (\bibinfo{year}{2019}) \bibinfo{pages}{1--14}.
\bibitem[{Kalhori et~al.(2017)Kalhori, Alamdari, Zhu, Samali, and
  Mustapha}]{kalhori}
\bibinfo{author}{H.~Kalhori}, \bibinfo{author}{M.~M. Alamdari},
  \bibinfo{author}{X.~Zhu}, \bibinfo{author}{B.~Samali},
  \bibinfo{author}{S.~Mustapha},
\newblock \bibinfo{title}{Non-intrusive schemes for speed and axle
  identification in bridge-weigh-in-motion systems},
\newblock \bibinfo{journal}{Measurement Science \& Technology}
  \bibinfo{volume}{28} (\bibinfo{year}{2017}) \bibinfo{pages}{025102}.
\bibitem[{Alamdari et~al.(2019)Alamdari, Kildashti, Samali, and
  Goudarzi}]{ex37}
\bibinfo{author}{M.~M. Alamdari}, \bibinfo{author}{K.~Kildashti},
  \bibinfo{author}{B.~Samali}, \bibinfo{author}{H.~V. Goudarzi},
\newblock \bibinfo{title}{Damage diagnosis in bridge structures using rotation
  influence line: Validation on a cable-stayed bridge},
\newblock \bibinfo{journal}{Engineering Structures} \bibinfo{volume}{185}
  (\bibinfo{year}{2019}) \bibinfo{pages}{1--14}.
\bibitem[{Besselink et~al.(2013)Besselink, Tabak, Lutowska, van~de Wouw,
  Nijmeijer, Rixen, Hochstenbach, and Schilders}]{ex34}
\bibinfo{author}{B.~Besselink}, \bibinfo{author}{U.~Tabak},
  \bibinfo{author}{A.~Lutowska}, \bibinfo{author}{N.~van~de Wouw},
  \bibinfo{author}{H.~Nijmeijer}, \bibinfo{author}{D.~Rixen},
  \bibinfo{author}{M.~Hochstenbach}, \bibinfo{author}{W.~Schilders},
\newblock \bibinfo{title}{A comparison of model reduction techniques from
  structural dynamics, numerical mathematics and systems and control},
\newblock \bibinfo{journal}{Journal of sound and vibration}
  \bibinfo{volume}{332} (\bibinfo{year}{2013}) \bibinfo{pages}{4403--4422}.
\bibitem[{Dormand and Prince(1980)}]{RK}
\bibinfo{author}{J.~Dormand}, \bibinfo{author}{P.~Prince},
\newblock \bibinfo{title}{A family of embedded runge-kutta formulae},
\newblock \bibinfo{journal}{Journal of Computational and Applied Mathematics}
  \bibinfo{volume}{6} (\bibinfo{year}{1980}) \bibinfo{pages}{19--26}.
\bibitem[{Kanungo et~al.(2002)Kanungo, Mount, Netanyahu, Piatko, Silverman, and
  Wu}]{ex38}
\bibinfo{author}{T.~Kanungo}, \bibinfo{author}{D.~M. Mount},
  \bibinfo{author}{N.~S. Netanyahu}, \bibinfo{author}{C.~D. Piatko},
  \bibinfo{author}{R.~Silverman}, \bibinfo{author}{A.~Y. Wu},
\newblock \bibinfo{title}{An efficient k-means clustering algorithm: analysis
  and implementation},
\newblock \bibinfo{journal}{IEEE Transactions on Pattern Analysis \& Machine
  Intelligence} \bibinfo{volume}{24} (\bibinfo{year}{2002})
  \bibinfo{pages}{881--892}.
\bibitem[{Yuan and Yang(2019)}]{ex39}
\bibinfo{author}{C.~Yuan}, \bibinfo{author}{H.~Yang},
\newblock \bibinfo{title}{Research on k-value selection method of k-means
  clustering algorithm},
\newblock \bibinfo{journal}{J} \bibinfo{volume}{2} (\bibinfo{year}{2019})
  \bibinfo{pages}{226--235}.
\bibitem[{Erturk and Inman(2009)}]{4.1}
\bibinfo{author}{A.~Erturk}, \bibinfo{author}{D.~J. Inman},
\newblock \bibinfo{title}{An experimentally validated bimorph cantilever model
  for piezoelectric energy harvesting from base excitations},
\newblock \bibinfo{journal}{Smart Materials \& Structures} \bibinfo{volume}{18}
  (\bibinfo{year}{2009}) \bibinfo{pages}{25009--25018}.
\bibitem[{Junior et~al.(2009)Junior, Erturk, and Inman}]{85}
\bibinfo{author}{C.~Junior}, \bibinfo{author}{A.~Erturk},
  \bibinfo{author}{D.~J. Inman},
\newblock \bibinfo{title}{An electromechanical finite element model for
  piezoelectric energy harvester plates},
\newblock \bibinfo{journal}{Journal of Sound and Vibration}
  \bibinfo{volume}{327} (\bibinfo{year}{2009}) \bibinfo{pages}{9--25}.

\end{thebibliography}
